%% file: HR_ieee_proc14_submitted_final_arxiv_revised.tex
\documentclass[12pt]{article}
\usepackage{amsmath}
\usepackage{bm}
\usepackage{graphicx}
\usepackage{color}
\usepackage[caption=false]{subfig}
\usepackage{multirow}
\usepackage{lscape}
\usepackage{amsfonts}
\usepackage{bm}
\usepackage{parskip}
\input{header_defns.tex}

\usepackage[tmargin=1.0in,bmargin=1.0in,rmargin=1.0in,lmargin=1.0in]{geometry}

\newtheorem{thm}{Theorem}
\newcommand{\con}{\mbox{\scriptsize con}}
\DeclareMathOperator*{\sign}{sign}

\begin{document}
\setlength{\parindent}{15pt}

\renewcommand{\thefootnote}{$\text{c1}$}

\begin{center}
\textbf{\large Foundational principles for large scale inference: Illustrations through correlation mining}\\
\end{center}
\vspace*{0.25in}
\begin{center}
Alfred O. Hero$^{\dagger}$ and Bala Rajaratnam$^{*}$
\end{center}
\vspace*{0.25in}

\noindent $^\dagger$ University of Michigan, Ann Arbor, MI 48109-2122, USA

\noindent $^*$ Stanford University, Stanford, CA 94305-4065, USA

\section*{Abstract}

When can reliable inference be drawn in the ``Big Data" context? This paper presents a framework for answering this fundamental question in the context of correlation mining, with implications for general large scale inference.  In large scale data applications like genomics, connectomics, and eco-informatics the dataset is often variable-rich but sample-starved: a regime where the number $n$ of acquired samples (statistical replicates) is far fewer than the number $p$ of observed variables (genes, neurons, voxels, or chemical constituents).  Much of recent work has focused on understanding the computational complexity of proposed methods for  ``Big Data". Sample complexity however has received relatively less attention, especially in the setting when the sample size $n$ is fixed, and the dimension $p$ grows without bound. To address this gap, we develop a unified statistical framework that explicitly quantifies the sample complexity of various inferential tasks. Sampling regimes can be divided into several categories: 1) the classical asymptotic regime where the variable dimension is fixed and the sample size goes to infinity; 2) the mixed asymptotic regime where both variable dimension and sample size go to infinity at comparable rates; 3) the purely high dimensional asymptotic regime where the variable dimension goes to infinity and the sample size is fixed. Each regime has its niche but only the latter regime applies to exa-scale data dimension.  We illustrate this high dimensional framework for the problem of correlation mining, where it is the matrix of pairwise and partial correlations among the variables that are of interest. Correlation mining arises in numerous applications and subsumes the regression context as a special case. We demonstrate various regimes of correlation mining based on the unifying perspective of high dimensional learning rates and sample complexity for different structured covariance models and different inference tasks. 

{\bf Keywords:} large scale inference, Big Data, sample complexity, asymptotic regimes, purely high dimensional, unifying learning theory, triple asymptotic framework, correlation mining, correlation estimation, correlation selection, correlation screening, graphical models.

\section{Introduction}
%
%
%
%
%
%
%
%
%
%
%
%
%
%
%
%

The increasing availability of large scale and high dimensional data is driving a major resurgence of data science, recently rebranded under the moniker ``Big Data" \cite{lynch2008big}. There has been a preponderance of catch phrases such as ``big-data-biology," ``ecoinformatics," ``precision medicine," ``data-driven decision-making," ``Big Data business analytics" in scientific publications and the media. However, until recently, most of the research in Big Data has concentrated on issues of data management, data warehousing,  computational data analysis, and end-user data utilization \cite{trelles2011big}, \cite{marx2013biology}, \cite{michener2012ecoinformatics}. While the data management research community has made progress on the problem of quality assurance, e.g., associated with provenance and computer errors \cite{labrinidis2012challenges}, the issue of limited sample size and statistical reproducibility remains largely open.  This issue has been recognized as one of the principal hurdles that stand in the way of success of the scientific enterprise \cite{young2011deming}, \cite{rudin2014discovery}.  The negative consequences of insufficient samples can be especially dire when the data is high dimensional, heterogenous and uncalibrated  \cite{dal2014making}, \cite{madigan2013evaluating}.  It is therefore both important and timely to address the problem of inference on Big Data from the point of view of statistical reproducibility.

The statistical reproducibility point of view is founded on a non-monolithic notion of Big Data: the data should be considered as a matrix with $p$ columns and $n$ rows indexed by, respectively, $p$ variables over each of the data fields and $n$ independent samples of these variables.  If there is only a single sample then the matrix collapses to a vector and no statistical analysis of reproducibility can be performed. With larger sample size, the probability of reproducibility can be studied in the context of the statistical theory of random matrices.

Here we develop this perspective for a particular Big Data problem: correlation mining in high dimension where the number of samples is much smaller than the dimension, a setting that we call ``sampled starved." Correlation mining is an area of data mining where the objective is to discover patterns of correlation between a large number ($p$) of observed variables based on a limited number ($n$) of samples.  Correlation mining can be framed as the mathematical problem of reliably reconstructing different attributes of the correlation matrix of the population from the sample covariance matrix that is empirically constructed from the $p \times n$ data matrix. The inference task depends on the attributes of the correlation that are of interest, while the performance of a correlation mining algorithm for a particular task depends on the number of samples and the underlying structure of the population covariance.  In high dimensional settings where $p \gg n$,  correlation mining presents significant challenges to the practitioner, both in terms of unavailability of computationally tractable algorithms and in terms of lack theory that could be used to specify sample size requirements. This paper will provide some perspective on the latter challenges.

Covariance matrices arise in a very diverse set of Big Data applications including: empirical finance and econometrics 
\cite{ledoit2003improved}, \cite{ledoit2004honey}, \cite{Ledoit&Wolf:JMA04}, \cite{todros2012measure}, \cite{bollen1998structural}, \cite{fornell1982two};   MIMO radar and communications \cite{guerci2000optimal}, \cite{gini2002covariance}, \cite{kim2001comparison}, \cite{bliss2002eim},
\cite{fuhrmann2008transmit}, \cite{li2009mimo},\cite{hero2003secure},  \cite{biglieri2007mimo};
image analysis \cite{geman&geman:1984},\cite{mardia1988multi},\cite{chellappa1985classification},\cite{Willsky:IEEE02};  network sensing \cite{vardi1996network}, \cite{chamberland2003decentralized}, \cite{patwari2005locating}, \cite{patwari2003relative}, \cite{chen2011robust}, \cite{chung2007detection}; life and biomedical sciences \cite{mcintosh1996spatial}, \cite{van1997localization}, \cite{almasy1998multipoint}, \cite{schafer2005shrinkage}, \cite{Korrat:2007ec},
\cite{zhu2005high}, \cite{rajapakse2010networking};  and climate science \cite{GRE2014}, \cite{WEGSR2014}, \cite{Schneider01}, \cite{Smerdon_et_al2010}, to name just a few. However, covariance matrices are used differently depending on the application and the task. For tracking of targets from space-time-adaptive-radar (STAP) the task is to estimate the entire covariance matrix in so far as it yields estimates of eigenvectors that span the signal (target) subspace \cite{guerci2003space}.  For exploring functional gene regulation networks the task is to determine the matrix locations of the largest pairwise correlations or inverse correlations, often obtained by thresholding the sample correlation matrix \cite{hero2011large} or by performing sparsity constrained optimization \cite{banerjee2006convex}, \cite{Friedman&etal:08}.  In linear discriminant analysis the task is to estimate a quadratic form of the inverse covariance matrix \cite{chen2011robust}, \cite{Morrison:90}. In anomaly detection, it is the Schur complement of the covariance matrix that is of interest as it is related to residual prediction error covariance \cite{Wiesel&Hero:SP09}. In independence screening one is interested in the support of the non-zero rows of the covariance \cite{fan2008sure}. In variable selection for prediction it is the support of the regression coefficient vector that is of interest \cite{Firouzi&etal:AISTATS13}
(a setting in which correlation mining subsumes the regression context as a special case).

Of central importance for all of these applications are the sample size requirements, which can differ from task to task. Reliably performing some of these correlation mining tasks might require relatively few samples, e.g., screening for the presence of variables that are hubs of high correlation in a sparsely correlated population,  while other tasks might require many more samples, e.g., accurately estimating all entries of the inverse covariance matrix in a densely correlated population. A theoretical framework has been emerging for predicting these sampling requirements as a function of the population correlation structure and as a function of the inference task. The principal aim of this paper is to present building blocks of this framework with an emphasis on the high dimensional setting where the number $p$ of variables is much larger than the number $n$ of samples ($p \gg n$), a setting relevant to correlation mining in massive data sets.

Several real-world examples of this high dimensional setting are given below:

\begin{itemize}
\item In studies of correlation networks of the annotated human genome $p$ is on the order of tens of thousands of genes while $n$ is typically fewer than a hundred samples, e.g., corresponding to a population of human subjects. In these studies correlation levels of magnitude as low as  $0.2$ are sometimes considered significant \cite{de2004discovery}, though some of them may actually be spurious.
\item In space-time-adaptive-processing (STAP) radar a spatio-temporal covariance matrix is used to filter out clutter in order to better detect a moving target at a particular range and Doppler frequency. For full degrees-of-freedom (DOF) STAP the estimator of the spatio-temporal  clutter covariance can have dimension $p=rq$ on the order of hundreds of thousands and the number of samples $n$ under a hundred. Here $r$ is the number of radar pulses (time), $q$ is the number of elements in the radar array (space), and $n$ is the number of range bins in the vicinity of the target \cite{guerci2003space},\cite{li2009mimo}.
\item For spambot network discovery from honeypot data, correlation of spamming profiles is performed between $p=$ hundreds of thousands of known IP addresses while the number of time points in each profile may only be on the order of $n=30$ \cite{Xu&etal:ICC2009}.
\item In recommender systems preference vectors of $p=$ millions of subscribers are correlated with other subscribers on the basis of $n=$ a few hundred preference categories, e.g., movies or music, to predict future user preferences \cite{koren2009matrix}, \cite{hsiao2014social}.
\item In fMRI connectomics a long term goal is to correlate brain activations across individual brain neurons, i.e.,  $p=10^{11}$.  Currently connectome researchers might parcel the brain  into $p$ regions, with $p=$ several thousand, and estimate correlation by averaging over $n=$  hundreds of repeated activation stimuli. The objective is to use the sample correlation matrix to study patterns of activation in order to reveal the ``brain network"  \cite{sporns2005human}, \cite{van2010exploring}.
\end{itemize}

Therefore, it is not an exaggeration to say that correlation mining practitioners face a deluge of variables with very limited sample size. With so few samples one is bound to find spurious correlations between some pairs of the many variables. It is therefore essential to understand the intrinsic sampling requirements of such statistical inference problems.  The study of sampling requirements falls into several different asymptotic regimes.  Classical statistical error prediction and control methods are based on an asymptotic regime where the dimension $p$ is fixed and the sample size $n$ goes to infinity. Such a regime is obviously inapplicable to the case of $n<p$. Recently, statisticians have developed theory that applies to the regime where $p$ goes to infinity. This theory is in the realm of high dimensional statistical inference, often called  the ``large $p$ small $n$ regime" \cite{west2003bayesian}, \cite{Buhlmann:2011}. This theory covers the case where both $n$ and $p$ go to infinity, which, as contrasted to the classical fixed $p$ regime, is a setting that we called the mixed asymptotic regime and
includes the so-called ``high dimensional," ``very-high dimensional" or ``ultra-high
dimensional" settings \cite{Buhlmann:2011}, \cite{fan2008sure}. Each of these denote regimes where the speed at which $n$ goes to infinity as a function  of $p$ becomes progressively slower.
However, since the theory requires that both $p$ and $n$ go to infinity it may not be very useful in applications where the availability of samples is limited and finite. In such a case a more useful and relevant regime is the ``purely-high dimensional" setting where the number of samples $n$  remains fixed while the number of variables $p$ goes to infinity. This setting applies, for example, to the correlation screening application \cite{hero2011large}, discussed in Sec. \ref{sec:screening}, where the correlation threshold $\rho$ approaches one as $p$ becomes large. This regime is in fact the highest possible dimensional regime, short of having no samples at all. Thus it is appropriate to call this purely-high dimensional regime the ``ultimately-high dimensional regime," and we shall use these two terms interchangeably in the sequel.
A table comparing these asymptotic regimes is given below
(Table \ref{table_asy}).

\begin{table}
\begin{center}
{\tiny
  \begin{tabular}{|c|c|c|c|c|c|}
    \hline
    \textbf{Asymptotic}& \textbf{Terminology} & \textbf{Sample}  &
\textbf{Model} & \textbf{Application} & \textbf{References} \\
    \textbf{framework} &\textbf{(``setting")} & \textbf{size} &
\textbf{dimension} & \textbf{domain} & \textbf{(selected)}\\
    \hline
    & & $n$ & $p$ & & \\
    \hline
    \multirow{6}{*}{
    \shortstack[c]{Large \\ sample \\ asymptotics
\\(Classical)}
    } & \multirow{6}{*}{\shortstack[c]{small \\ dimensional}} &
\multirow{6}{*}{$\longrightarrow \infty$} & \multirow{6}{*}{fixed} &
\multirow{6}{*}{``small data''} & Fisher \cite{Fisher1922, Fisher1925},
Rao \cite{Rao1947, Rao1963}, \\
    & & & & & Neyman and Pearson \cite{NeymanPearson1933}, Wilks
\cite{Wilks1938}, \\
    & & & & & Wald \cite{Wald1941a, Wald1941b, Wald1943, Wald1949}, \\
    & & & & & Cram\'{e}r \cite{Cramer1946a, Cramer1946b}, Le Cam
\cite{LeCam1953, LeCam1986}, \\
    & & & & & Chernoff \cite{Chernoff1956}, Kiefer and
Wolfowitz\cite{KieferWolfowitz1956},  \\
    & & & & & Bahadur \cite{Bahadur1967}, Efron \cite{Efron1982} \\
    \hline
    \multirow{9}{*}{\shortstack[c]{Mixed \\asymptotics}} &
\multirow{9}{*}{\shortstack[c]{high \\ to \\ ultra high \\dimensional}} &
\multirow{9}{*}{$\longrightarrow \infty$} &
\multirow{9}{*}{$\longrightarrow \infty$} &
\multirow{9}{*}{\shortstack[c]{``medium sized'' data \\ (mega or giga
scales)} } & \multirow{9}{*}{\shortstack[c]{Donoho \cite{Donoho2006},
Zhao and Yu \cite{ZhaoYu2006},  \\ Meinshausen and B\"{u}hlmann
\cite{Meinshausen2006a}, \\
Cand\`{e}s and Tao \cite{CandesTao2007}, Wainwright
\cite{Wainwright2009a, Wainwright2009b},\\
Bickel, Ritov, and Tsybakov\cite{BickelRitovTsybakov2009},\\
Peng, Wang, Zhou, and Zhu \cite{peng2009a},  \\
Khare, Oh, and Rajaratnam, \cite{KhareOhRaj2014},\\ Fan and Lv \cite{fan2008sure} }} \\
  & & & & & \\
  & & & & & \\
  & & & & & \\
  & & & & & \\
  & & & & & \\
  & & & & & \\
  & & & & & \\
  & & & & & \\
    \hline
    \multirow{4}{*}{\shortstack[c]{Large \\ dimension \\
asymptotics}}
    & \multirow{4}{*}{\shortstack[c]{purely high \\dimensional}} &
\multirow{4}{*}{fixed} & \multirow{4}{*}{$\longrightarrow \infty$} &
\multirow{4}{*}{\shortstack[c]{Sample Starved \\ ``Big Data'' \\ (tera to
exascales)} } &  \multirow{4}{*}{\shortstack[c]{Hero and Rajaratnam \cite{hero2011large}, \\ Hero and Rajaratnam \cite{hero2012hub}, \\ Firouzi, Hero and Rajaratnam
\cite{FirouziHeroRajaratnam2014} }}\\
  & & & & & \\
  & & & & & \\
  & & & & & \\
   \hline
  \end{tabular}
  }
\end{center}
\caption{\em \small Overview of different asymptotic regimes of
statistical data analysis. These regimes are determined by the relation
between the number $n$ of samples drawn from the population and the
number $p$, called the dimension, of variables. In the classical
asymptotic regime the number $p$ is fixed while $n$ goes to infinity.
This is the regime where most of the well known classical statistical
testing procedures, such as student $t$ tests of the mean, Fisher $F$
tests of the variance, and Pearson tests of the correlation, can be
applied reliably. Mixed asymptotic regimes where $n$ and $p$ both go to
infinity have received much attention in modern statistics. However, in
this era of big data  where $p$ is exceedingly large, the mixed
asymptotic regime is inadequate since it still requires that $n$ go to
infinity.  The recently introduced ``purely high dimensional regime"
\cite{hero2011large}, \cite{hero2012hub} 
addressed in this paper,
is more suitable to big data problems where $n$ is limited and finite
(adapted from \cite{Hero&Rajaratnam:Cambridge14})}.
\label{table_asy}
\end{table}


The purely-high dimensional regime of large $p$ and fixed $n$ is a mathematical characterization of the extremely big data problem. Evidently, this regime poses several challenges to correlation mining practitioners. These include both  computational challenges and the challenge of error control and performance prediction. Yet this regime also holds some rather pleasant surprises. Remarkably, there  are surprising benefits to having few samples in terms of computation and scaling laws.  There is a scalable computational complexity advantage relative to other high dimensional regimes where both $n$ and $p$ are large. In particular,  correlation mining algorithms can take advantage of numerical linear algebra shortcuts and approximate $k$ nearest neighbor search to compute large sparse correlation or partial correlation networks. Another benefit of purely-high dimensionality is an advantageous scaling law of the false positive rates as a function
of $n$ and $p$. Even small increases in sample size can provide significant gains in this  regime. For example, when the dimension is $p=10,000$ and the number of samples is
$n=100$, the experimenter only needs to double the number of samples in order to
accommodate an increase in dimension by six orders of magnitude ($p=10,000,000,000$)
without increasing the false positive rate \cite{hero2012hub}.

The unifying comparative tool used in this paper is sample complexity analysis. We will use this analysis to place different statistical models and different inference tasks into the various asymptotic regimes shown Table \ref{table_asy} (see Tables \ref{table:contextual}-\ref{table:tasks}). Sample complexity analysis has been widely applied to study the performance of inference procedures.  In supervised computational learning sample complexity is defined as the number of training samples necessary to maintain a given level of generalization error as a function of the complexity of the statistical model or algorithm \cite{bartlett1998sample} \cite{haussler1994bounds}, \cite{niyogi1996relationship}, \cite{dasgupta2005coarse}, \cite{kakade2003sample}.  Similar notions arise in the context of high dimensional convergence rates
\cite{Buhlmann:2011}, Bayesian model complexity \cite{spiegelhalter2002bayesian}, and stochastic complexity (minimum description length - MDL) \cite{rissanen1989stochastic}\footnote{Stochastic complexity is different from the notion of statistical complexity in statistical mechanics \cite{PhysRevLett.63.105}}. There is always a tradeoff between statistical sample complexity and statistical model complexity. The nature of this tradeoff is specific to the model and to the task. This specificity forms the basis for the unifying the theory of large scale inference that is developed in this paper.

The outline of the paper is as follows. Sec. \ref{sec:correlation} reviews the definitions of covariance, precision, correlation, partial correlation, and Gaussian graphical models. Section \ref{sec:estimation} treats the problem of covariance estimation. Section \ref{sec:model_selection} develops the problem of model selection and support estimation. In Sec. \ref{sec:screening} the problem of correlation screening is discussed. In  Sec. \ref{sec:complexity} sample complexity comparisons of various types of correlation inference tasks are presented. Concluding remarks are given in Sec. \ref{sec:conclusions}.

%
\section{Correlation matrices and correlation networks}
\label{sec:correlation}

\subsection{Covariance and precision matrices}
\label{sec:covprecision}

Let $\bX=[x_1, \ldots, x_p]^T \in \Reals^p$ be a real-valued random (column) vector having a distribution  whose second  order moments exist. The mean vector $\bmu=E[\bX] \in \Reals^p$ and the covariance matrix $\bSigma=\cov(\bX)=E[(\bX-\bmu)(\bX-\bmu)^T] \in \Reals^{p\times p}$ are defined in terms of the statistical expectation operator $E[\cdot]$ associated with the distribution of $\bX$. The matrix $\bSigma$ is symmetric positive semidefinite and its $ij$'th entry is the covariance $\cov(X_i,X_j)= E[(X_i-\mu_i)(X_j-\mu_j)]$, with $\mu_i$ the $i$-th element of $\bmu$. It will be assumed throughout that $\bSigma$ is in fact positive definite, which is generally true when the $p$ components of $\bX$ are non-redundant random variables, i.e., no variable can be predicted without error using linear combinations of the other variables. The inverse of the covariance matrix is called the precision matrix  $\bOmega= \bSigma^{-1}$. The covariance and precision matrices capture marginal dependency and conditional dependency, respectively,  in the joint distribution of $\bX$. In particular, when $\bSigma$  ($\bOmega$) is sparse, i.e., it has only a few non-zero elements, there are many variables that are marginally (conditionally) independent of the other variables.  More will be said about this in Sec. \ref{sec:GMM}.

The  matrices $\bSigma$ and $\bOmega$ are not invariant to scaling of the component variables of $\bX$. Scale invariance is important when one is interested in a measure of dependency that is not a function of the units used to represent different variables. The correlation matrix $\bR$ and partial correlation matrix $\bP$ are respective versions of $\bSigma$ and $\bOmega$ that are scale invariant
\begin{eqnarray}
\bR&=& \diag(\bSigma)^{-1/2} \bSigma \diag(\bSigma)^{-1/2}
\label{eq:Rdef} \\
\bP&=& \diag(\bOmega)^{-1/2} \bOmega \diag(\bOmega)^{-1/2}
\label{eq:Pdef}
\end{eqnarray}
where, for a $p\times p$ matrix $\bB$, $\diag(\bB)$ is the $p\times p$ diagonal matrix formed from the diagonal elements of $\bB$ and, for a diagonal matrix $\bD$ with non-zero valued entries $\{d_{ii}\}_{i=1}^p$ along the diagonal,  $\bD^{-1/2}$ denotes the diagonal matrix with entries $\{1/\sqrt{d_{ii}}\}_{i=1}^p$ along the diagonal. The diagonal elements of $\bR$ and $\bP$ are equal to $1$ and the off-diagonal elements are between $-1$ and $1$. An important property of $\bR$ and $\bP$ is that they retain any zero patterning in $\bSigma$ and $\bOmega$, respectively,  if $\bSigma$ or $\bOmega$ are sparse.

Let $\{\bX_k\}_{k=1}^n$ be a set of $n$ i.i.d. realizations from the distribution of $\bX$.
The sample covariance matrix is the $p \times p$ positive semidefinite matrix
\begin{equation}
\bS_n= \frac{1}{n-1} \sum_{k=0}^n (\bX_k- \hat{\bmu})(\bX_k-\hat{\bmu})^T
\label{eq:Sn}
\end{equation}
where $\bmu=n^{-1} \sum_{k=0}^n \bX_k$ is the sample mean vector. The sample covariance $\bS_n$ is a statistically consistent estimator of $\bSigma$, i.e., it converges to $\bSigma$ with probability one as $n\rightarrow\infinity$ if $p$ is fixed. Furthermore, the plug-in estimator $\hat{\bR}=\diag(\bS_n)^{-1/2} \bS_n \diag(\bS_n)^{-1/2}$ can be used to estimate $\bR$. The matrix $\bS_n$ is positive definite with probability one if $n\geq p$ and if the distribution of $\bX$ is Lebesgue continuous. In this case, plug-in estimators $\hat{\bOmega}=\bS_n^{-1}$ and $\hat{\bP}=\diag(\hat{\bOmega})^{-1/2}\hat{\bOmega} \diag(\hat{\bOmega})^{-1/2}$ can be used to consistently estimate $\bOmega$ and $\bP$. In the case $n<p$, which is the main focus of this paper, alternative estimation strategies must be adopted.

\subsection{The Gaussian graphical model}
\label{sec:GMM}
Assume that the data is a set of $n$ i.i.d. realizations $\{\bX_k\}_{k=1}^n$ of a random vector $\bX\in \Reals^p$ that has a multivariate Gaussian distribution
with positive definite covariance matrix $\bSigma$.  When the precision matrix  $\bOmega$ is sparse such data is said to be a Gauss Markov random field (GMRF) or, equivalently, is said to be distributed according to a Gaussian graphical model (GGM).
The log-likelihood function for $\bSigma$ and $\bOmega$ can be expressed as
\be
l(\bSigma)&=&- \frac{n}{2} \trace\{\bSigma^{-1} \bS_n\} -\frac{n}{2} \logdet (\bSigma),\\
\label{eq:GMMll}
l(\bOmega)  &=& -\frac{n}{2}\trace\{\bOmega \bS_n\} +\frac{n}{2} \logdet (\bOmega),
\label{eq:GMMlp}
\ee
up to unimportant additive constants. In the GGM $\bOmega$ is a more natural parameterization than $\bSigma$ in the sense that it is sparse while $\bSigma$ may not be sparse. In addition, the Gaussian Graphical model corresponds to a natural exponential family with $\bOmega$ as the canonical parameter. Furthermore, the sparse structure of $\bOmega$ is directly related to conditional independencies specified by the pairwise, local and global Markov properties of undirected graphical models and also by the factorization of the joint density of $\bX$ \cite{lauritzen1996graphical},\cite{wainwright2008graphical}. 

The GGM is a ``graphical model" in the following sense. Define a graph $\mathcal G$
whose $p$ vertices (or nodes) correspond to the $p$ variables $X_1, \ldots, X_p$ in $\bX$. Assign an edge between node $i$ and node $j$ if and only if the $ij$-th entry of the sparse precision matrix $\bOmega$  is non-zero.  The edges of the so-constructed graph $\mathcal G$, often called a correlation or partial correlation graph, are specified by an adjacency matrix $\bA$. $\bA$  is a binary matrix that indicates the {\em support} of the matrix $\bOmega$; specifically the non-zero entries of $\bA$ indicate the locations the non-zero entries of $\bOmega$.  The degree of node $i$ is the number of neighboring nodes in the graph and corresponds to the number of non-zero entries in the $i$-th row of $\bOmega$. In
GMRFs it is typically assumed that the graph $\mathcal G$ has only a small number
($O(p)$) of edges, corresponding to the property $\bOmega$ is row sparse\footnote{A $p\times p$ matrix is row sparse with sparsity coefficient $s$ if no row has more than $s$ non-zero entries, where $s=o(p)$.} \cite{lauritzen1996graphical}.

Observe that in a GGM it is the inverse covariance $\bOmega$ and not the covariance itself that is important. There are many other reasons for the wide interest in $\bOmega$ and its scaled version $\bP$, the partial correlation matrix (\ref{eq:Pdef}) \cite{Muirhead:82}, \cite{Anderson:03}, \cite{lauritzen1996graphical}, \cite{Buhlmann:2011},\cite{VanTrees}, \cite{Kay:98},  \cite{Kay:91}. First, given predictor variables $\bX=[X_1, \ldots, X_n]^T$, the optimal linear predictor of a response variable $Y$  depends on $\cov(\bX)$ only through its inverse.  Second, when $\bX$ is Gaussian,  the Bayes-optimal test between two hypothesized covariances similarly depends only on the inverse covariances. Third, many classical multivariate procedures such as linear discriminant analysis (LDA) and multivariate analysis of variance (MANOVA) depend directly on the inverse covariance matrix.  Fourth, the $ij$-th entry of $\bP$ is equal to the conditional correlation between $X_i$ and $X_j$ given the rest of the variables $\{X_k\}_{k\neq i,j}$. Fifth, if $\hat{X}_i$ is the optimal linear predictor of $X_i$ given $\{X_k\}_{k\neq i}$, the $ij$-th entry of the partial correlation matrix $\bP$ is equal to the correlation coefficient between the prediction error residuals associated with $\hat{X}_i$ and $\hat{X}_j$.  Finally, many physical random processes have a sparse inverse covariance matrix while the covariance matrix is not sparse. We illustrate this last point by the following example.

\subsection{Illustrative example of a GGM}

We illustrate the differences between the covariance and inverse covariance matrices
in a GGM using an example motivated by statistical physics. Let
$X(u_1,u_2)$ be a function on the unit square  $[0,1]^2 \subset \Reals^2$ that satisfies
the  Poisson partial differential equation
$$\begin{array}{cc}
\Delta X(u_1,u_2)=W(u_1,u_2), & u_1, u_2 \in  (0,1)^2 \\
X(u_1,u_2)= 0, & u_1, u_2 \not \in (0,1)^2 \end{array}$$
where $W(u_1,u_2)$ is a fixed integrable function. Here $\Delta X=\frac{\partial^2
X}{\partial u_1^2} + \frac{\partial^2 X}{\partial u_2^2}$ is the Laplacian differential
operator. Solutions $X$ to the Poisson equation arise in heat transfer, electromagnetics,
and fluid dynamics \cite{sommerfeld1949partial}.

We extract a Gaussian graphical model from the Poisson equation by discretizing it to a
finite difference equation and setting the discretized $W$ to be a Gaussian random
matrix.  This will convert the continuous domain function $X$ over $[0,1]^2$ to a vector
$\bX$ with $p=N_1 N_2$ elements, where $N_1=1/\delta_1$ and $N_2=1/\delta_2$ and
$\delta_1$, $\delta_2$ are the $u_1$ and $u_2$ increments used in the finite difference
approximation. Specifically, discretize $u_1$ and $u_2$ to the sets $u_1 \in
\{k\delta_1\}_{k=0}^{N_1}$ and $u_2 \in \{k\delta_2\}_{k=0}^{N_2}$. Then, assuming smooth
$W$,  the discretized $X$ satisfies the finite difference equation:
\begin{equation}
X_{i,j}=\frac{(X_{i+1,j}+X_{i-1,j}) \delta^2_2+(X_{i,j+1}+X_{i,j-1}) \delta^2_1
-W_{i,j}\delta_1 \delta_2}{2(\delta^2_1+\delta^2_2)}.
\label{eq:dPoisson}
\end{equation}
Second, arrange the $\{X_{i,j}\}_{i,j=1}^{N_1,N_2}$ into the vector
$\bX\in \Reals^{N^2}$
in lexicographic order.  This results in an equivalent matrix form of the equation
(\ref{eq:dPoisson}) $$[\bI-\bA]\bX=\bW$$ where $\bA$ is a sparse tri-diagonal matrix.

\begin{figure}[h!]
        \centering
       \includegraphics[width=6.3cm]{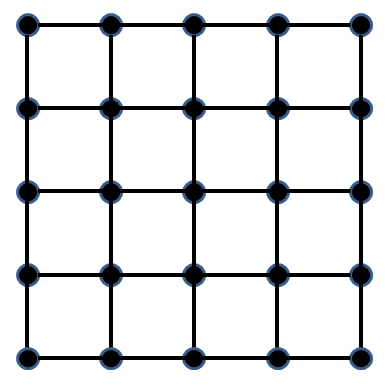} \hspace{0.2in}
       \includegraphics[width=9cm]{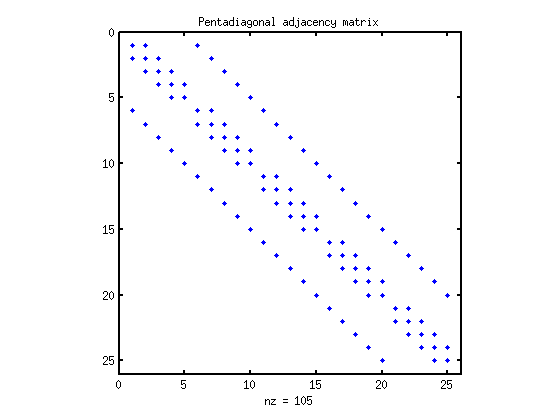}
       \caption{\em \small The Gaussian graphical model obtained from a finite difference
       approximation to the Poisson equation has only local spatial dependency (left
       panel) resulting in a pentadiagonal inverse covariance matrix (right panel).}
\label{fig:randomfield}
\end{figure}

Under the assumption that  $\bW$ is a i.i.d. zero mean spatio temporal Gaussian driving noise, i.e., $\cov(\bW)= \sigma_W^2 \bI$,
the covariance $\bSigma$ and the inverse covariance
$\bOmega$ of $\bX$ have the form $$\bSigma = \bOmega^{-1}, \hspace{0.2in} \bOmega=
[\bI-\bA][\bI-\bA]^{T}/\sigma_W^2.$$
As $\bA$ is sparse
tri-diagonal, the inverse covariance $\bOmega$ and partial correlation $\bP$ are both sparse penta-diagonal matrices with support set shown in Fig. \ref{fig:randomfield} for $N_1=N_2=5$. Therefore,  the random vector $\bX$ is a Gauss Markov random field with GGM spatial dependency structure specified by the non-zero entries of $\bOmega$, equivalently of $\bP$. Note that the covariance matrix $\bSigma$ is
not sparse as the inverse of the penta-diagonal matrix  is a full matrix.  Thus,  as is usually the case for a GGM, the inverse covariance is a more parsimonious model description than is the covariance.

\begin{figure}
\begin{center}
\includegraphics[height=5cm]{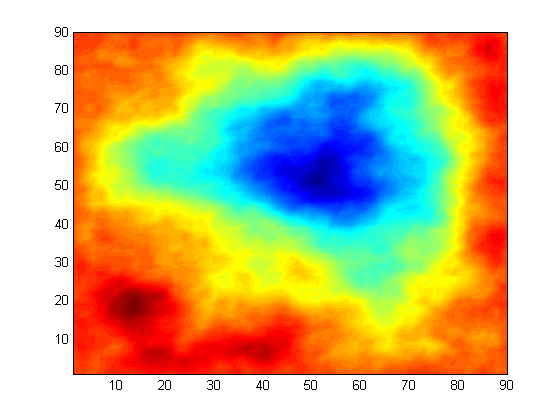} \;
\includegraphics[height=5cm]{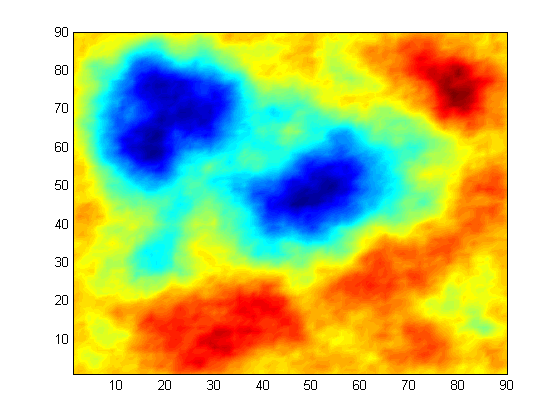}
\end{center}
\caption{\em \small Two ($n=2$) independent realizations of the discretized Poisson random field obtained from the finite difference approximation (\ref{eq:dPoisson}) with $N_1=N_2=90$. Here the driving process $\bW$ is a zero mean spatio-temporally white Gaussian noise.
}
\label{fig:realizationPoisson}
\end{figure}

\begin{figure}
\begin{center}
\includegraphics[height=12cm]{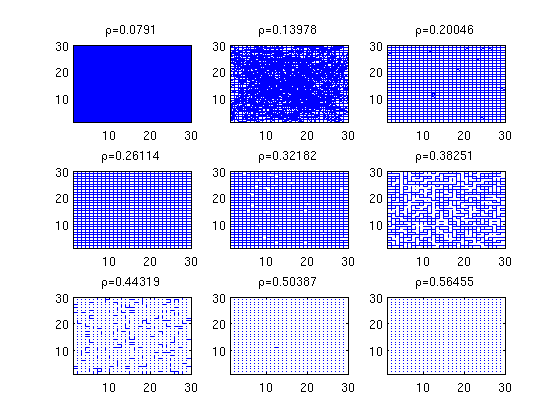}
\end{center}
\caption{\em \small The support of the thresholded sample partial correlation matrix, rendered as a graph over spatial locations in the plane, for a Gaussian random field generated by the Poisson equation on a $30\times 30$ grid ($p=900$). The sample partial  correlation $\hat{\bP}=\diag(\hat{\bOmega})^{-1/2} \hat{\bOmega} \diag(\hat{\bOmega})^{-1/2}$, where $\hat{\bOmega}=\bS_n^{-1}$ the inverse of the sample covariance matrix, was computed using $n=1500$ samples.   The graph is shown for six different threshold levels $\rho$ applied to $\hat{\bP}$.  The true cartesian support set is recovered for $\rho=0.26114$. There is a sharp increase in the number of false edges as the threshold $\rho$ decreases below a certain threshold, somewhere between 0.0791 and 0.13978.  The location of this threshold appears well approximated by the theory \cite{hero2011large} that predicts a critical threshold value $0.091$ (See equation (\ref{eq:rhocrit}) in Sec. \ref{sec:screening} of this paper).  }
\label{fig:parcorest}
\end{figure}

For visualization, two realizations of $\bX$ are shown in Fig.  \ref{fig:realizationPoisson} for a simulation of the discretized Poisson equation (\ref{eq:dPoisson}) on a $30 \times 30$ spatial grid ($p=900$). While general methods of support set estimation are discussed in Sec. \ref{sec:model_selection}, Fig. \ref{fig:parcorest} illustrates a simple empirical estimator of the support set of $\bOmega$ based on $n=1500$ samples of $\bX$.   This simple estimator is constructed by thresholding an estimate $\hat{\bP}$ of the partial correlation matrix $\bP$ at a level $\rho \in[0,1]$, which is user defined. The non-zero entries of the adjacency matrix associated with the thresholded $\hat{\bP}$ specifies the support set estimator. Here $\hat{\bP}=\diag(\hat{\bOmega})^{-1/2} \hat{\bOmega}\diag(\hat{\bOmega})^{-1/2}$, where $\hat{\bOmega}=\hat{\bS}^{-1}_n$ is the inverse of the sample covariance $\bS_n$ defined in (\ref{eq:Sn}).   Figure  \ref{fig:parcorest} shows the support set estimator for six different values of the applied threshold $\rho$. Note that when $\rho=0.26114$ the true support set, illustrated in Fig. \ref{fig:randomfield} for a $5\times 5$ grid, is recovered. When $\rho$ decreases below a critical value there is an abrupt increase in the number of false edges. Below this critical value (top left panel of Fig. \ref{fig:parcorest}) the number of false edges   approaches the number ${p \choose 2}$ of edges in the complete graph.

%
%
%
%
%

\section{Correlation mining: estimation}
\label{sec:estimation}
%
%
%

An ambitious correlation mining objective is to empirically estimate each and every one of the entries of the covariance $\bSigma$ or inverse covariance matrix $\bOmega$ from the $n$ samples. The accuracy of these estimators is often measured by the
Frobenius norm of the estimation error. Other common criteria for gauging estimator accuracy are the matrix $\ell_2$ norm,  also called the operator norm,   and the  matrix $\ell_1$ norm. For example, in  estimation of $\bOmega$ the Frobenius norm difference between $\bOmega$ and its empirical estimate $\hat{\bOmega}$ is:
$$\|\bOmega-\hat{\bOmega}\|_F=\sqrt{\sum_{i,j=1}^p (\omega_{ij}-\hat{\omega}_{ij})^2},$$
where $\omega_{ij}$ denotes the $ij$-th entry of $\bOmega$.
When the objective is to estimate the
entries of the covariance, the correlation or the partial correlation matrix, the Frobenius norm error is defined similarly.   Covariance and inverse covariance estimation arise in many applications including:
finance \cite{Ledoit&Wolf:JMA04}, \cite{ghaoui2003worst};  gene expression modeling
\cite{zhu2007bayesian}, \cite{rothman2008sparse}; sensor array processing
\cite{Hero&Delap:Haykin95},\cite{Kelly&Forsythe:TR848}, \cite{Johnson&Dudgeon:93};
space-time-adaptive processing (STAP) radar \cite{wicks2006space}, \cite{guerci2003space}
\cite{stoica2008using}; spatio-temporal classification \cite{Greenewald2013} and
prediction \cite{Tsiligkaridis&Hero:TSP2013Kroneckerdecomposition};  brain connectomics
\cite{craddock2013imaging},\cite{van2012human}; and sensor network anomaly detection
\cite{chen2011robust}.

When the data are multivariate Gaussian,  the  maximum likelihood (ML)
estimator of $\bOmega$ maximizes the log-likelihood function
(\ref{eq:GMMlp}).
In the case when $n\geq p$ the ML estimator of the covariance $\bSigma$ is equal to $\hat{\bSigma}=\frac{n-1}{n}\bS_n$, where $S_n$ is the sample covariance matrix defined in (\ref{eq:Sn}), and the ML estimator of $\bOmega$ is $\hat{\bOmega}=\hat{\bSigma}^{-1}$. In high dimensional situations where $n <p$, $\bS_n$
is not positive definite and these ML estimators do not exist. Estimator regularization
is commonly used  in  this situation in order to yield a positive definite solution.
Model-based regularization methods add a suitable penalty function $c(\bOmega)$
to the log-likelihood function resulting in the so-called
penalized ML estimator.  When the penalty is interpreted as a log prior on $\bOmega$ the
penalized ML estimator is equivalent to the maximum a posteriori (MAP) estimator,
also known as the posterior mode estimator.  Commonly used regularization
penalties are the quadratic penalty $\lambda\|\bOmega\|_F^2$ and the non-sparsity penalty
$\lambda \|\bOmega\|_1=\lambda \sum_{i\neq j}^p |\omega_{ij}|$, where $\lambda>0$ is called
the regularization parameter, and pushes the penalized ML solution towards a diagonal matrix
as $\lambda$ increases to $\infty$. Other penalties shrink the penalized ML solution towards more highly
structured matrices, e.g.,  block sparse, Toeplitz, banded, circulant, or Kronecker matrix structure.
We call the latter  structured penalties as contrasted to the un-structured penalties that encourage
the estimator to be diagonal. When the true inverse covariance $\bOmega$ has the same structure as the structure induced by the penalty
we say that the penalty of the penalized ML estimator is matched to the model.

The effect of different structured models for $\bOmega$ on  the asymptotic convergence
rate of the penalized ML error can be  of interest to practitioners.   In large scale
settings the high dimensional convergence rate is a natural comparative performance measure. High dimensional rates of convergence have been obtained under
a large range of structured and unstructured penalties (and matched models).
Often, when the penalty is matched to the model, the Frobenius norm estimation error decreases at asymptotic
rate of the form $\sqrt{P(\log p)/n}$ as $p,n \rightarrow \infty$ where $P$ is the
(effective) number of free parameters in the model.  We illustrate for the case of spatio-temporal multivariate Gaussian graphical models; a poster child for high dimensional inference.

\subsection{Illustration: estimation of a spatio-temporal precision matrix}

Consider the case where a set of $q$ sensors acquires $r$ time samples (a snapshot) of a spatio-temporal random field, e.g., the STAP radar example \cite{guerci2003space},\cite{li2009mimo} discussed in the introduction.   There are $n$ snapshots of this random field. Each snapshot generates a multivariate Gaussian distributed random matrix $\mathbb X$ with $q$ rows (spatial) and $r$ columns (temporal).  Define the lexicographically ordered vectorization  $\bX={\mathrm{vec}}(\mathbb X) \in \Reals^{qr}$ of $\mathbb X$. The  spatio-temporal covariance  $\cov(\bX)=\bSigma$ of a snapshot is  an unknown $p \times p$ matrix where $p=qr$. The task is to estimate the precision matrix $\bOmega=\bSigma^{-1}$.
When there is no additional information about $\bSigma$, beyond that it is symmetric and positive definite,  the precision matrix $\bOmega$ is completely unknown and there are  $P={p \choose 2}=O(p^2)$ unknown parameters to be estimated, which is the number of distinct entries of the matrix. The base case of no information  is often called the saturated model. If the experimenter is given additional contextual information  about the structure of the precision matrix the value of this information can often be quantified in terms of sample complexity. We give several examples below.

\noindent{\bf Spatio-temporal contextual information: sparse precision matrix}:
Contextual information specifying that $\mathbb X$ is a GGM with sparse inverse covariance matrix results in a significant reduction in the number of free parameters in $\bOmega$:  from $P=O(p^2)$ to $P=O(p)$. The problem of sparse covariance estimation has been of significant interest \cite{Meinshausen2006a}, \cite{Banerjee2008a}, \cite{Friedman&etal:08}. When the spatio-temporal matrix $\mathbb X$ represents the outputs of sensor network, an example of
contextual information is the knowledge that the physical environment is a time-varying random field whose amplitudes obey the laws of Lagrangian classical mechanics, e.g., fluid flow, heat flow, or electro-magnetic wave fields that satisfy Poisson or Navier-Stokes partial differential equations \cite{Schafer&etal:PhysicsofScales}. The contextual information may simply be an upper bound on the sparsity parameter $s$ or it
may actually specify the sparsity pattern, i.e., fix the graph $\mathcal G$ specifying the support of $\bOmega$.

\noindent{\bf Spatio-temporal contextual information: Kronecker covariance}:
Any spatio-temporal covariance has some degree of Kronecker structure that can be
captured by the Kronecker sum decomposition of Pitsianis and Van Loan
\cite{van1993approximation}, \cite{Tsiligkaridis&Hero:TSP2013Kroneckerdecomposition}. When applied to the inverse covariance, this decomposition is $\bOmega=\sum_{k=1}^k \bA_i \otimes \bB_i$ where $k$ is the Kronecker rank and $\bA_i, \bB_i$ are linearly independent $q\times q$ and $r\times r$ matrices, called Kronecker factors. When the contextual information is  that $\bOmega$ has Kronecker rank $k$ the number of free parameters in $\bOmega$ is reduced from $O(p^2)$ to $O(k(q^2+r^2))$.  The simplest case occurs when $\bOmega$ is known to have Kronecker rank 1: $\bOmega=\bA \kron \bB$, where $\bA$ and $\bB$ are symmetric positive definite
$q\times q$ (spatial) and $r\times r$ (temporal) covariance matrices. Covariance estimation in this Kronecker model, also known as the matrix normal model, has been widely studied \cite{Dutilleul1999}, \cite{werner2008estimation},  \cite{tsiligkaridis2013convergenceTSP}. A physical example of Kronecker-structured covariance is the case that $\mathbb X$ are outputs of a passive
antenna array that is sensing the EM emissions of $k$ wide-sense-stationary non-moving
targets in the far field of the array.  In particular, if $k=1$ there is a
single target and the spatial and temporal components of the covariance are completely
decoupled. In this case, the contextual information could be the number of targets.

\noindent{\bf Spatio-temporal contextual information: Kronecker+sparse precision matrix}: Consider the case that the contextual information specifies that $\mathbb X$ has a precision matrix $\bOmega$ that is both Kronecker structured and has sparse  Kronecker factors. When the Kronecker rank is $k=1$ then $\bOmega$ factors into $\bA\otimes\bB$ and there are $P=O(qs_A+rs_B)$ unknown parameters, where
$s_A$ and $s_B$ specify the sparsity in each of the Kronecker factors $\bA$ and $\bB$
\cite{tsiligkaridis2013convergenceTSP}, \cite{Allen2010}.

The value of the different types of contextual information can be studied in terms of
estimation of $\bOmega$ using theory developed recently in
\cite{tsiligkaridis2013convergenceTSP}, which discusses estimators and estimator convergence rates for the cases that $\bOmega$ is, respectively, sparse, Kronecker, and
sparse-Kronecker structured. These structural constraints are incorporated into the
estimator $\hat{\bOmega}$  of $\bOmega$ using penalized maximum likelihood methods, which can be interpreted as MAP covariance estimators when the prior
distributions on the entries of $\bOmega$ are specified by the penalty functions.   For
each  case denote the mean-squared error of the MAP estimator as the
expectation of the Frobenius norm difference squared:
MSE$=E[\|\bOmega-\hat{\bOmega}\|_F^2]$. The relative decrease in the associated MSE due
to any of the above pieces of contextual information is a function of the sample complexity for the inverse covariance estimation task. The forms of these MAP estimators, their asymptotic log MSE, and their asymptotic sample complexity can be obtained directly from \cite{tsiligkaridis2013convergenceTSP},
\cite{Tsiligkaridis&Hero:TSP2013Kroneckerdecomposition} and are summarized below and in
Table \ref{table:contextual}.

For the saturated model (no side information about $\bOmega$) the maximum likelihood
estimator is equal to the inverse of the sample covariance $\hat\bSigma= n^{-1}
\sum_{i=1}^n (\bx_i-\overline{\bx})(\bx_i-\overline{\bx})^T$, where $\overline{\bx}$ is
the sample mean of  $\{\bx_i\}_{i=1}^n$ \cite{Morrison:90}. This is the MAP
estimator under the uniform estimator-loss function and a uninformative (constant) prior on $\bOmega$. When contextual information specifies that $\bOmega$ is in fact sparse, the Glasso penalized ML estimator  \cite{Friedman&etal:08} adds an l1-norm penalty on $\bOmega$ to the log-likelihood function. The Glasso precision estimator is the MAP estimator under a Laplacian-like prior on $\bOmega$ and it  can be determined using an iterative maximization algorithm. When contextual information is that $\bOmega$ has Kronecker structure, again an iterative algorithm must be used to find the maximum
likelihood estimator under the Kronecker model \cite{Werner2007}. This is the Bayes
optimal estimator under a Kronecker covariance model and non-informative priors on the
Kronecker factors.   Finally, when the contextual information is that $\bOmega$ is both
Kronecker and sparse the Kronecker log-likelihood model can be penalized by l1-norms on
each of the Kronecker factors, resulting in a MAP estimator under a
Laplacian-type  prior on each of the factors \cite{Allen2010},
\cite{tsiligkaridis2013convergenceTSP}.

The value-added brought to estimator performance by each of these contextual information sources can be assessed by studying the asymptotic sample complexity. The
asymptotic sample complexity is defined by the number $n=n_p$ of samples, as a function
of the dimension $p$, required to maintain a given level of performance as $p$ goes to
infinity. Under general conditions on the MAP spatio-temporal covariance
estimators, the log MSE takes the high dimensional (large $q,r$ and large $n$) asymptotic form 
shown in the 3rd row of Table \ref{table:contextual} \cite{tsiligkaridis2013convergenceTSP}. The second row of the table gives the form for
the prior on $\bOmega$ as (from left to right): 1) uniform prior in the case of no
contextual information; 2) sparse $\ell_1$ prior in the case of sparsity information on
$\bOmega$; 3) Kronecker structured $\bOmega$ in the case of contextual information that
decouples space and time (the prior $\delta(\rank {\mathrm R}(\bOmega)-1) $ is a delta
function that forces kronecker structure $\bOmega=\bA\otimes\bB$)\footnote{$\mathcal R$
is the permutation-rearrangement operator
\cite{Tsiligkaridis&Hero:TSP2013Kroneckerdecomposition} that maps $qr\times qr$ matrix
$\bOmega$ into the $q^2\times r^2$ matrix ${\mathcal R}(\bOmega)$, which has rank 1 if
and only if $\bOmega =\bA\otimes \bB$ for some $q \times q$ matrix $\bA$ and $r \times r$
matrix $\bB$.}; 4) Kronecker plus sparse $\bOmega$.

\begin{table}
{\small
\centering
\begin{tabular}{|c||c|c|c|c|}
   \hline
   {\bf {Information}}   & None & Sparse & Kronecker & Kronecker+sparse \\
   \hline \hline
  {\bf {Model}} &  saturated $\bOmega$&  sparse $\bOmega$&  $\bOmega =\bA\otimes \bB$ & sparse $\bOmega=\bA\otimes \bB$ \\
  \hline
    {\bf {log}} $f(\bOmega)$ & constant & $\lambda \|\bOmega\|_1$  & $\delta\left(\rank{\mathcal R}(\bOmega)-1\right)$  & $\delta\left(\rank{\mathcal R}(\bOmega)-1\right)
    +\lambda_1 \|\bA\|_1+\lambda_2 \|\bB\|_2$  \\
   \hline
   {\bf Bound} & $\frac{1}{2}\log\left(\frac{q^2r^2}{n}\right)$ & $\frac{1}{2}\log\left(\frac{qr\log qr}{n}\right)$ & $\frac{1}{2}\log\left(\frac{(q^2+r^2)\log M }{n}\right)$  & $\frac{1}{2}\log\left(\frac{(q+r)\log M}{n}\right)$  \\
   \hline
   {\bf {Regime}} & $\frac{q^2r^2}{n}\rightarrow \alpha$ & $\frac{qr\log qr}{n}\rightarrow\alpha$ & $\frac{(q^2+r^2)\log M}{n}\rightarrow\alpha$  & $\frac{(q+r)\log M}{n}\rightarrow \alpha$  \\
   \hline
\end{tabular}
}
\caption{\em \small
Expressions  for performance of MAP estimators (4th row) and the associated regimes of asymptotic sample complexity (5th row) for estimating the $qr \times qr$ spatio-temporal inverse correlation matrix $\bOmega =\bSigma^{-1}=(\cov(\bX))^{-1}$, associated with the $q \times r$ space-time Gaussian random matrix $\mathbb X$,  for different types of models of $\bOmega$ representing prior contextual information (2nd and 3rd rows). The contextual information specifies patterns of sparsity and/or dependency between the $q$ spatial coordinates and the $r$ temporal  coordinates. The bound is the asymptotic (large $q,r$ and $n$) log Frobenius norm error of the Bayes-optimal MAP estimator of $\bOmega$, for the different types of contextual information (1st row). The four rightmost columns of the table correspond respectively to: no contextual information about $\bOmega$ (None); information that $\bOmega$ is sparse, corresponding to $\mathbb X$ being a Gauss Markov random field (GMRF), with $\lambda>0$ sufficiently large so that the sparsity factor is of order $o(qr)$; information that $\bOmega$ has (rank $k=1$) Kronecker product structure (Kronecker), corresponding to decoupled rows and columns of $\mathbb X$; and information that $\bOmega$ has (rank $k=1$) Kronecker product structure with sparse Kronecker factors (Kronecker GMRF), with $\lambda_1,\lambda_2 >0$ ufficicently large so that sparsity factors are of order $o(max\{q,r\})$.   Comparing  the  sample complexity regimes (4th row) between the various types of contextual information quantifies the value of taking additional samples (See Fig. \ref{fig:kronVoI}).}
\label{table:contextual}
\end{table}

\begin{figure}[h!]
\centering
\includegraphics[width=8cm]{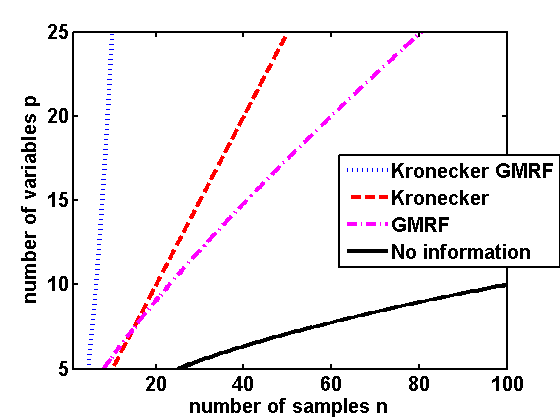}
\includegraphics[width=8cm]{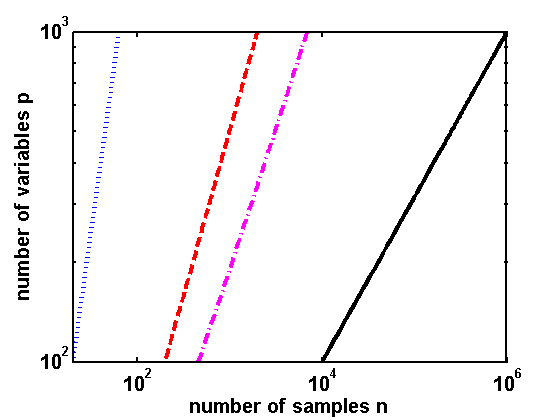}
\caption{\em \small Sample complexity for estimating the $qr
\times qr$ spatio-temporal correlation matrix for different types of prior contextual
information shown in the 1st row of Table \ref{table:contextual} for the case
$q=r=\sqrt{p}$, where $\sqrt{p}$ is a positive integer. The curves show asymptotic sample
complexity (5th row in the Table), which are constant contours of the proxies (4th row of the
Table) over the plane of the number $p$ of variables and the number of samples $n$. The
asymptotic proxies are equal over all the curves and along each curve. Curves to the left
represent lower sample complexity and indicate the reduction in the required
number of samples to attain a given level of log MSE for specified number of parameters.
Here the contextual information represents knowledge that the inverse covariance has
sparse structure alone (curve labeled ``GMRF"), Kronecker structure alone (curve labeled
``Kronecker"), vs Kronecker+GMRF structure. The curve labeled ``No information"
represents the case where there is no contextual side information about the inverse
covariance. The curve labeled Kronecker GMRF dominates all of the others since
information that the inverse covariance has both Kronecker and GMRF structure achieves
maximal reduction in the number of free parameters and provides the highest value per
sample.}
\label{fig:kronVoI}
\end{figure}

By quantifying the change in log MSE associated with different types of contextual
information (sparse, Kronecker, Kronecker+sparse) the value of taking additional samples can be determined across the contextual information regimes shown in Table \ref{table:contextual}. For the purpose of comparison we assume that $q=r=\sqrt{p}$, i.e., the spatial and temporal dimensions are identical. For the different categories of contextual information we fix the number of variables $p$ and the log MSE. Figure \ref{fig:kronVoI} plots the level sets of constant log MSE over the number of samples and the number of variables. These curves indicate that the knowledge of both sparse and Kronecker structure is more valuable than knowledge of either sparse or Kroenecker structure alone. To illustrate, assume that $\mathbb X$ is a $100 \times 10$ matrix corresponding to $10$ shapshots of $100$ sensors.  Then $p=1000$ and, from the right panel of Fig. \ref{fig:kronVoI}, if the contextual information specified that the inverse covariance has Kronecker and sparse structure the MAP estimator requires only $n=75$ samples as compared to $4000$ or $8000$ samples if Kronecker or sparse structure is specified, and $n=1,000,000$ samples if no information were available. The value of the information that the covariance is both sparse and Kronecker structure is that it decreases sampling requirements by more than $4$ orders of magnitude relative to no contextual information!

\section{Correlation mining: model selection}
\label{sec:model_selection}

One of the primary goals of model selection in the correlation mining setting is to identify the support of the correlation or partial correlation matrix, i.e., identify the pairs of variables with non-zero correlations or partial correlations, from $n$ measurements of the set of $p$ variables \cite{Banerjee2008a}, \cite{Friedman&etal:08}, \cite{Hsieh2011}, \cite{Guillot2012}, \cite{DalalRaj2014}, \cite{Meinshausen2006a}, \cite{peng2009a}, \cite{Rocha2008}, \cite{Friedman2010}, \cite{LeeHastie2014}, \cite{KhareOhRaj2014}, \cite{ODKRNips14}. Model selection should be easier than estimating the values of all the correlations, discussed in Sec. \ref{sec:estimation}.  The original \emph{covariance selection} problem \cite{Dempster1972} considers estimating inverse covariance matrices with zero entries under a multivariate Gaussian model for the observations. In recent years, the problem of estimating the support of sparse inverse covariance matrices has become a popular topic in the high dimensional statistics and machine learning literature.

As discussed in Sec. \ref{sec:covprecision}, the problem of identifying the
sparse patterning structure in the inverse covariance matrix $\Omega$ is equivalent to identifying the graph $\mathcal G$ associated with the non-zero entries of $\Omega$,  and is thus also popularly known as graphical model selection. Various approaches have been proposed for identifying graphical models from high dimensional data. They can be  categorized broadly into penalized likelihood methods and Bayesian methods. Popular Bayesian methods entail specifying priors on the space of sparse covariance or inverse covariance matrices \cite{dawid_lauritzen},  \cite{letac_massam_2007}, \cite{Rajaratnam2008}, \cite{Khare2011} and using Bayesian scoring rules to undertake model selection. Penalized likelihood methods in the context of (partial) correlation mining can be further categorized into a) Gaussian-based penalized likelihood methods \cite{Banerjee2008a}, \cite{Friedman&etal:08}, \cite{Hsieh2011}, \cite{Guillot2012}, \cite{DalalRaj2014}, and b) Pseudo-likelihood based or regression based methods \cite{Meinshausen2006a}, \cite{peng2009a}, \cite{Rocha2008}, \cite{Friedman2010}, \cite{LeeHastie2014}, \cite{KhareOhRaj2014}, \cite{ODKRNips14}.

Recent work on penalized likelihood methods have focused on understanding both the computational complexity and sample complexity of model selection approaches. The former entails computing the sparse inverse covariance estimate using $\ell_1$-penalized likelihood approaches and identifying the corresponding graph. This includes among others developing fast iterative algorithms for maximizing $\ell_1$-penalized likelihoods in order to obtain sparse inverse covariance estimates, quantifying the computational complexity of these algorithms, and deriving convergence rates of the iterative algorithms. In line with the main theme of this paper, the focus of this section will be to understand the sample complexity of model selection problems in the correlation mining context.

Sample complexity of model selection approaches are generally stated in terms of sign consistency and estimation consistency (see \cite{peng2009a}, \cite{KhareOhRaj2014}). As in covariance estimation, the set-up is to let both the sample size $n$, and the dimension $p=p_n$ tend to infinity and establish large sample properties for a sequence of covariance parameters that is growing in dimension. We shall illustrate the sample complexity of model selection via a concrete example below. In particular, we present below a sample complexity result of a recently proposed graphical model selection method proposed in \cite{KhareOhRaj2014} called CONCORD.

CONCORD seeks to maximize the following jointly convex objective function, called the pseudo-likelihood, as proposed in \cite{KhareOhRaj2014}:
\begin{align}
  Q_{\con} (\Omega) &=: -\sum_{i=1}^p n\log \omega_{ii} + \frac{1}{2} \sum_{i=1}^p \|
  \omega_{ii} {\bf Y}_i + \sum_{j \neq i} \omega_{ij} {\bf Y}_j \|_2^2
  + \lambda \sum_{1 \leq i < j \leq p} |\omega_{ij}|, \label{eq2}
\end{align}
where $\omega_{ij}$ is the $ij$-th element of the $p \times p$ matrix $\Omega$, and $\mathbf Y_i$ denotes the $i$-th feature vector.
Iterative optimization algorithms, along with aspects of computational complexity and algorithm convergence, are covered in \cite{KhareOhRaj2014}, \cite{KRopt2014},  \cite{ODKRNips14}. For sample complexity results, we shall follow very closely large sample results of the CONCORD graphical model selection approach as given in \cite{KhareOhRaj2014}. Both estimation consistency and oracle properties under suitable regularity conditions are stated below. The reader is referred to \cite{KhareOhRaj2014} for further technical details.

Let $\{\bar{\Omega}_n\}_{n \geq 1}$ denote the sequence of true underlying inverse covariance matrices and let the dimension $p = p_n$ vary with the sample size $n$. As in \cite{peng2009a}, assume the existence of accurate estimates of the diagonal entries $\{\widehat{\alpha}_{n,ii}\}_{1
\leq i \leq p_n}$ such that for any $\eta > 0$, there exists a constant $C > 0$ such that
$$
\max_{1 \leq i \leq p_n} \left| \widehat{\alpha}_{n,ii} - \bar{\omega}_{ii} \right| \leq C \left(
\sqrt{\frac{\log n}{n}} \right),
$$
holds with probability larger than $1 - O(n^{-\eta})$.

For vectors $\omega^o \in \mathbb{R}^{\frac{p_n(p_n-1)}{2}}$ and
$\omega^d \in \mathbb{R}^{p_n}_+$, let the notation ${\mathcal{L}}_n
(\omega^o, \omega^d)$ denote $\frac{{\mathcal{L}}_{\con}}{n}$
evaluated at a matrix with off-diagonal entries $\omega^o$ and diagonal entries $\omega^d$.
Let $\bar{\omega}_n^o =((\bar{\omega}_{n,ij}))_{1 \leq i < j \leq
  p_n}$ denote the vector of off-diagonal entries of $\bar{\Omega}_n$,
and $\widehat{\boldsymbol \alpha}_{p_n} \in \mathbb{R}^{p_n}_+$
denotes the vector with entries $\{\widehat{\alpha}_{n,ii}\}_{1 \leq i
  \leq p_n}$.  Let the sequence $\mathcal{A}_n$ denote the set of non-zero entries
in the vector $\bar{\omega}_n^o$, and let $q_n = |\mathcal{A}_n|$. Let $\bar\theta_{n,ij} =
\frac{\bar\omega_{n,ij}}{\sqrt{\widehat{\alpha}_{n,ii} \widehat{\alpha}_{n,jj}}}$ for $1 \leq i < j \leq p_n$ and define $\bar{\theta}_n^o = ((\bar{\theta}_{n,ij}))_{1 \leq i < j \leq   p_n} \in \mathbb{R}^{{p_n}({p_n}-1)/2}$. Also, let $s_n =
\min_{(i,j) \in \mathcal{A}_n} \bar{\omega}_{n,ij}$.   Assume furthermore that the following three regularity conditions are met (i) the spectrum of  $\bar{\Omega}_n$ is uniformly bounded from above and below, (ii) sub-Gaussianity of the data, (iii) the incoherence condition \cite{Meinshausen2006a}.

Under the above assumptions the following model selection result is established in \cite{KhareOhRaj2014}.
\begin{thm}\cite{KhareOhRaj2014}
  Suppose that assumptions (i), (ii), (iii) are satisfied. Suppose $p_n =
  O(n^\kappa)$ for some $\kappa > 0$, $q_n = o \left( \sqrt{n/\log n}
  \right)$, $\sqrt{\frac{q_n \log n}{n}} = o(\lambda_n)$, $\lambda_n
  \sqrt{n/\log n} \rightarrow \infty$, $\frac{s_n}{\sqrt{q_n} \lambda_n}
  \rightarrow \infty$ and $\sqrt{q_n} \lambda_n\rightarrow 0$, as $n
  \rightarrow \infty$.  Then there exists a constant $C$ such that for
  any $\eta > 0$, the following events hold with probability at least
  $1 - O(n^{-\eta})$.
  \begin{itemize}
  \item There exists a minimizer $\widehat{\omega}^o_n =
    ((\widehat{\omega}_{n,ij}))_{1 \leq i < j \leq p_n}$ of the CONCORD objective function $Q_{\con}
    (\omega^o, \widehat{\boldsymbol \alpha}_n)$.
  \item Any minimizer $\widehat{\omega}^o_n$ of $Q_{\con} (\omega^o,
    \widehat{\boldsymbol \alpha}_n)$ satisfies $ \| \widehat{\omega}^o_n
    - \bar{\omega}^o_n \|_2 \leq C \sqrt{q_n} \lambda_n $ and $
    \sign(\widehat{\omega}_{n,ij}) = \sign(\bar{\omega}_{n,ij}), \;
    \forall\ 1 \leq i < j \leq p_n.$
  \end{itemize}
\label{KhareOhRaj2014_thm}
\end{thm}

\noindent

The above result in \cite{KhareOhRaj2014} establishes model selection consistency (or rather sign consistency to be precise), and is in spirit similar to other model selection consistency results in the literature (see also \cite{Meinshausen2006a} and \cite{peng2009a}). A few remarks are in order with regards to sample complexity. First, note that model selection consistency requires that both the sample $n$ and dimension $p_n$ tend to infinity.  Hence asymptotic guarantees require large sample sizes and thus model selection consistency may not be valid in sample starved settings. Second, results in model selection consistency are often proved under the assumption of sub-Gaussianity of the tails, and thus may be restrictive in many applications with heavy-tailed data. Third, note that the dimension $p_n$ can grow faster than the sample size $n$, but cannot grow faster than a polynomial rate.

\section{Correlation mining: screening}
\label{sec:screening}

%
%
%

In correlation screening  one seeks to discover patterns of high correlation or partial correlation
between $p$ variables based on a set of $n$ observations \cite{hero2011large},
\cite{hero2012hub}. Stated in terms of a Gaussian graphical model, the objective is to
infer topological characteristics of the graph $\mathcal G$ associated with the zeros in
the precision matrix $\bOmega$.  In \cite{hero2011large}, we treated the problem of screening
for the presence of variables with high correlations to other variables. In \cite{hero2012hub},
we considered the setting of screening for the presence of connected nodes and hubs in $\mathcal G$ with
high partial correlation. Screening for such topology characteristics of $\mathcal G$
should be easier than model selection or covariance estimation.
Similarly to what was demonstrated  for covariance estimation (recall Table \ref{table:contextual}),
contextual information can be of high value for screening. For example, one may be given information
that specifies a
certain sparsity level or block diagonal structure of the inverse covariance, or information on
the minimum level of correlation that exists among the active variables in a block.

The correlation screening method of 
\cite{hero2012hub} finds edges, hubs, and other subgraph structures of $\mathcal G$ by
performing hypothesis testing. The method applies a
threshold to an empirical estimate $\hat{\bP}$ of the partial correlation matrix $\bP$,
defined in (\ref{eq:Pdef}), placing an edge in $\mathcal G$ where the magnitude of the  entry of $\hat{\bP}$ exceeds the threshold.
When $n\geq p$,  $\hat{\bP}$ may be the simple plug-in estimator of the
matrix inverse of the sample correlation estimator, while if $n< p$ the correlation screening methods developed in \cite{hero2012hub} uses the Moore-Penrose
generalized inverse of the sample correlation estimator.  Correlation screening has been studied and applied to hub discovery \cite{hero2012hub}, edge discovery
\cite{Hero&Rajaratnam:Cambridge14} and classification of local node degree
\cite{firouzi2013local} in a variety of graphical model applications including:
stationary Gaussian spatio-temporal processes models \cite{Firouzi&etal:GlobalSIP13},
\cite{firouziExcursions2014}; sparse regression models \cite{Firouzi&etal:AISTATS13};
and multiple model testing for common sparsity patterns \cite{hero2011large}.

The computational complexity of correlation screening is much lower than model selection or covariance estimation,
only of the order of $O(n^2)$, in the sample starved case of $n\ll p$.
To illustrate, the sample complexity of screening edges in $\mathcal G$ can be determined from the following
theorem (adapted from \cite[Prop. 2]{hero2012hub}, see also \cite{Hero&Rajaratnam:Cambridge14}):

\begin{thm}\label{nodescreentheorem}
Assume that the $n$ samples $\{\bX_k\}_{k=1}^n$ are i.i.d.
random vectors in $\mathbb R^p$ with bounded elliptically contoured density and
block sparse $p \times p$ covariance matrix. Let $\rho$ be the threshold applied
to the sample partial correlation matrix $\hat{\bP}$. Assume that $p$ goes to $\infty$ and
$\rho=\rho_p$ goes to one at a rate specified by the relation
$\lim_{p\rightarrow \infty} p(p-1)(1-\rho_p^2)^{(n-2)/2} = e_{n}$, where $e_n\in(0,\infty)$ is given.
Then the probability $P_e$ that there exists at least one false edge in $\mathcal G$ satisfies
\begin{equation}
\lim_{p\rightarrow \infty}P_e=1-\exp(-\kappa_{n}/2),
\label{eq:Pe}
\end{equation}
where
$$\kappa_{n}= e_n a_n (n-2)^{-1}$$
where $a_n$ is the volume of the $n-2$ dimensional unit sphere in $\Reals^{n-1}$
and is given by $a_n=\frac{\Gamma((n-1)/2)}{\sqrt{\pi} \Gamma((n-2)/2)}$.
\end{thm}

As in other types of hypothesis testing problems, two types of correlation screening
errors can occur: false positives (Type I) and false negatives (Type II). It is common
for an experimenter to constrain the false positive rate to ensure a certain level of Type I error control. Remarkably, Thm. \ref{nodescreentheorem} asserts that, under the stated conditions, the large $p$ false positive rate  does not depend on the true covariance $\bSigma$. In this large $p$ case the correlation threshold $\rho$ can be set to attain a given level of false positive
control. This fortunate situation is analogous to constant false alarm rate (CFAR) signal
detection in radar processing  \cite{Schwartz:IT69},
\cite{Robey&etal:AES92},\cite{Reed&Yu:ASSP90}.  In \cite{hero2011large} it was
shown that, when  screening a large number $p$ of variables, the false positive rate
undergoes a fundamental phase transition as a function of the applied threshold: the rate precipitously increases from almost to zero to one as the correlation threshold is
decreased beyond a critical threshold $\rho_c$. A direct consequence of Thm. \ref{nodescreentheorem}
is that $\rho_c$ has the form \cite[Eq. (10)]{hero2012hub}:
\begin{equation} \rho_{c}=\sqrt{1-(a_{n}(p-1))^{-2/(n-4)}},
\label{eq:rhocrit}
\end{equation}
When the applied threshold is greater than $\rho_c$ there will be few
false positives while when it is below $\rho_c$ the system will be inundated by false
positives.

\begin{figure}[h!]
\centering
\includegraphics[width=8cm]{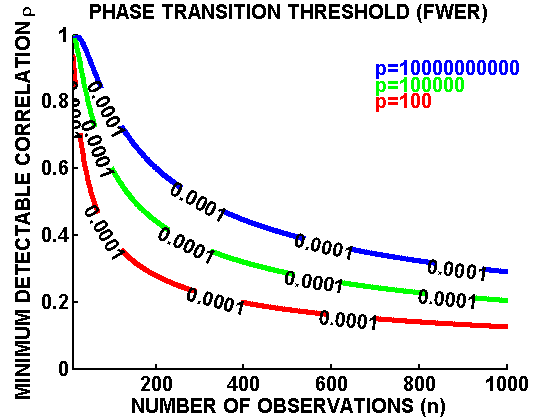}
\includegraphics[width=8cm]{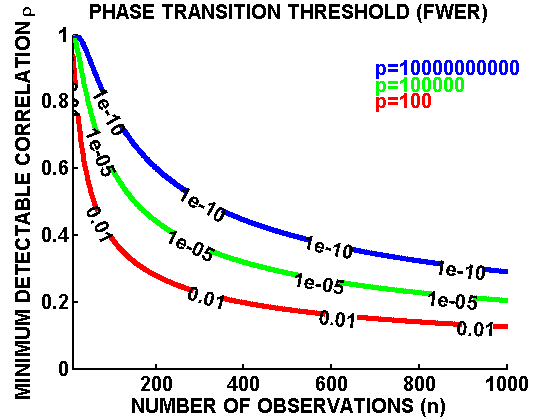}
\caption{\em \small Correlation screening curves quantifying value of information
associated with acquiring more samples $n$ for different parameter dimensions ($p=100,
10,000,$ and $10,000,000,000$) in terms of the minimum detectable correlation value
$\rho$.  The screening task is to detect variables that have high correlation (greater
than $\rho$) to at least one other variable. These curves specify the minimum required
number of samples $(n)$ for reliable detection of such variables at given family wise
false positive error rates. For example, for ten billion ($10^{10}$) variables at least
$200$ samples are required to reliably detect a variable having correlation greater than
$\rho=0.6$, while fewer than half the number of samples would be needed to detect the
same level of correlation if there were only ten thousand variables. Thus, for the
correlation screening task,  the value of a sample is much higher when there are fewer
variables, and the displayed curves quantify this value.  The curves in left panel are
isoclines on the probability of error surface for fixed family wise error rate (FWER)
equal to $0.0001$. The curves in the right panel are similar except that they are
isoclinal for fixed mean false positive rate of $1$ (only $1$ false positive node
detected out of $p$ nodes). The respective FWER's of false positive probability are
designated on each curve in the right panel.  The curves in the left and right panel are
very similar since the probability of error    surface undergoes an abrupt phase
transition from $0$ to $1$.}
\label{fig:rhovsn_phase_transition}
\end{figure}

\begin{figure}[h!]
\centering
\includegraphics[width=10cm]{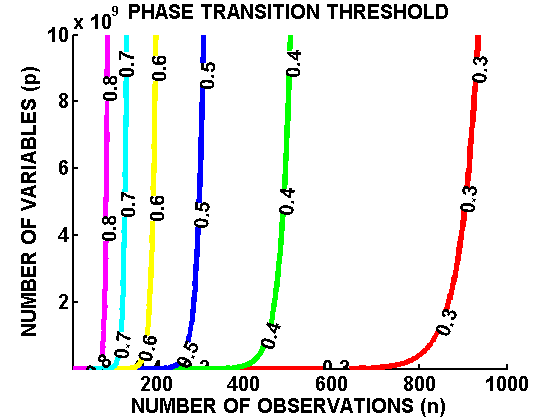}
\caption{\em \small Correlation screening curves quantifying value of information
associated with acquiring more samples $n$ for different minimal detectable correlation
levels ($\rho=0.3,0.4,0.5,0.6,0.7,0.8$) in terms of the parameter dimension $p$.  The
screening task is the same as in Fig. \ref{fig:rhovsn_phase_transition} but the phase
transition is plotted differently to reveal the value of information for detecting
different fixed levels of correlation for varying numbers of parameters $p$. Note that
the number of samples $n$ required for reliably detecting variables with high
correlations, e.g., $\rho=0.8$, increases much more slowly as $p$ increases than it does
for small correlations, e.g., $\rho=0.3$. Thus, as the desired correlation level
increases, there is a diminishing return in the value of information delivered by
acquiring additional samples.}
\label{fig:pvsn_phase_transition2}
\end{figure}

Using (\ref{eq:rhocrit}), and a large $p$ approximation \cite[Eq. (22)]{hero2012hub} to the false positive
rate  following directly from Thm. \ref{nodescreentheorem}, we can quantify the intrinsic value of taking an additional sample when the task is to  detect variables that have true correlations
exceeding a given threshold $\rho$. Assume that $p$ is fixed but large. Figures
\ref{fig:rhovsn_phase_transition} and \ref{fig:pvsn_phase_transition2} gives a family of
design curves that can be used by the system designer to right-size $n$ for given $p$ and given desired correlation level $\rho$.

\section{Correlation mining: intrinsic sample complexity regimes}
\label{sec:complexity}
%

The experimenter is often faced with several correlation mining tasks, possibly performed in sequence. For example, detection of existence of high correlations among $p$ variables may be followed by identification of the set of highly correlated variables, followed by estimation of the values of their correlations, followed by specification of the uncertainty (confidence intervals) associated with these estimates.
Reliable accomplishment of each task becomes  more difficult as one progresses from
detection to uncertainty quantification, requiring progressively larger numbers of samples and
progressively smaller critical phase transition thresholds $\rho_c$.
Establishing the sampling regimes associated with each one of these tasks is one of the fundamental
problems of large scale inference and data science.


Recall that the asymptotic sample complexity associated with an inference  task is the number $n=n_p$ of samples, as a function of the dimension $p$, required to maintain a given value of risk as $p$ goes to infinity.   Table \ref{table:tasks} summarizes
the sample complexity regimes (3rd row) and critical phase transition
threshold regimes (4th row) for tasks relevant to correlation mining. These are discussed below in order of increasing sample
complexity.
\begin{itemize}
\item {\bf Screening}: The screening task is to use the sample correlation to detect the existence of high correlations, by which we mean large values in either the population correlation or partial correlation matrix; equivalently to detect the existence of an edge in the correlation or partial correlation network, as discussed in Sec. \ref{sec:screening}. As such, it is a binary hypothesis testing problem having risk function equal to the false positive probability  under the null hypothesis $H_0$ of a sparse and invertible population covariance matrix $\bSigma$. For the threshold-based correlation screening method of \cite{hero2011large} the false positive probability is $P(N_e>0)$, where $N_e$ is the number of entries (edges) of the  sample partial correlation matrix that exceed a threshold $\rho$.   The bound is the asymptotic limiting value of the false positive edge probability specified in Thm. \ref{nodescreentheorem}.  
The sample complexity regime is:
$n$=fixed (not a function of $p$) while $p\rightarrow \infty$, denoted in Table
\ref{table:tasks} as $\frac{\log p}{n}\rightarrow \infty$. The critical threshold
(\ref{eq:rhocrit}) converges to $1$ as $p\rightarrow \infty$.
\item {\bf Detection}: The detection task is the same as the screening task except that both high and low correlations are of interest. The experimenter specifies a threshold $\rho\in (0,1)$ and the objective is to find correlations of magnitude at least $\rho$. The sample complexity regime for this problem is: $n$ increases to infinity with $p$ at the asymptotic rate $\frac{\log p}{n} =\alpha$, where $\alpha\in (0,\infty)$ is a constant. The critical threshold (\ref{eq:rhocrit}) converges to a constant $\rho^*$ satisfying $0<\rho^*<\rho$ as $p\rightarrow \infty$.
\item {\bf Support recovery}: This is the problem of model selection discussed in Sec. \ref{sec:model_selection} where the objective is to identify the support set $\calS\subset \{1, \ldots, p\}$ of the population inverse covariance $\bOmega$.
If $\hat{\calS}$ denotes an estimator of this support set, the risk function is the
probability that indices in $\calS$ are missing in $\hat{\calS}$ or that indices in $\{1, \ldots, p\}\subset \calS$ are erroneously included in $\hat{\calS}$, denoted $P({\mathrm {card}}\{{\calS} \Delta \hat{\calS}\}=\phi)$ where $A\Delta B$ is the symmetric difference between sets $A$ and $B$. Assume it is a priori known that the cardinality of $\calS$ is at most $k$, where $1\leq k \leq p$.
A finite sample upper bound follows by applying the union bound over the possible subsets of $\{1, \ldots, p\}$ of cardinality at most $k$ to obtain: $P({\mathrm {card}}\{{\calS} \Delta \hat{\calS}\}=\phi)\leq \sum_{i=0}^k {p \choose i} e^{-\beta n}$, where $\beta$ is the minimum Kullback-Liebler divergence between these subsets.
The well known bound \cite{arratia1989tutorial} on partial sums of binomial coefficients $\sum_{i=0}^k {p \choose i}\leq 2^{H(k/p)p}$, where $H(\epsilon)=-\epsilon \log \epsilon -(1-\epsilon) \log (1-\epsilon)$, can then be used along with the representation $H(k/p)p = p^\nu$ for some $\nu \in (0,1]$ to obtain the risk bound:  $2^{p^\nu}
e^{-\beta n} $.
Therefore, the limiting regime of values $(n,p)$ for which this bound is constant gives the sample complexity:  $n$ increases to infinity with $p$ at the asymptotic rate $\frac{p^\nu}{n} =\alpha$, where $\nu \in (0,1)$ and  $\alpha\in (0,\infty)$. Note that this is consistent with the rate $p=O(n^\kappa)$ of the CONCORD support recovery algorithm, given in  Thm. \ref{KhareOhRaj2014_thm}, for $\kappa>1$, the regime of where the number of samples is sufficient for model selection but not for parameter estimation. In this regime, the critical threshold (\ref{eq:rhocrit}) converges to zero as $p\rightarrow \infty$. Note that the asymptotic rates reported in Table \ref{table:tasks} for support recovery (or equivalently model selection) appear to be slower than what has been derived in the statistics literature. In particular, results are available which assert that provided $ (\log p )/n \to 0$, support recovery is possible with probability tending to 1 (see \cite{Rothman2009} for more details). We note that the these faster rates are a direct consequence of assuming that the variables are either Gaussian or  sub-Gaussian, or by imposing some tail condition \cite{Rothman2009}. When such conditions are relaxed, the asymptotic rates coincide with the regimes given in Table \ref{table:tasks} (see also \cite{elkaroui:annstat08, bickel&levina:aos08_2} for details on convergence in matrix norms when Gaussianity is relaxed).

\item {\bf Parameter estimation}: The problem of parameter estimation is to determine the individual values of the $p\choose 2$ entries in
$\bOmega$. The risk function is the MSE, defined as the mean Frobenius norm squared error $E[\|\bOmega-\hat\bOmega\|^2_F]$
between the population inverse covariance  and an empirical estimator. The bound is the high dimensional limiting value of this MSE as $n\rightarrow \infty$ and
$p\rightarrow \infty$ \cite{Buhlmann:2011}. The sample complexity  for this problem is: $n$ increases to infinity with $p$ at the asymptotic rate $\frac{p \log p}{n} =\alpha$,
where $\alpha\in (0,\infty)$ is a constant. Again the critical threshold
(\ref{eq:rhocrit}) converges to zero as $p\rightarrow \infty$.
\item {\bf Performance estimation}:  We consider the most general (and stringent) setting
for performance estimation where, for a specified Borel set $\mathcal B \subset \Reals^p$, the
probability $P(\bX\in \mathcal B)$ must be accurately estimated. For example, assuming
that $\bX$ has a zero mean elliptically contoured density $f(\bx)$, if $\mathcal B$ is
the set $\{\bx \in \Reals^p: \|\bx\| >\gamma\}$ for $\gamma >0$, $P(\bX\in \mathcal B)$
is the critical region for optimally rejecting outliers and detecting anomalies
\cite{Scott&Nowak:JMLR06}, \cite{Hero:NIPS06} and the value $\hat{\gamma}=\|\bX\|$  can
be used to estimate the p-value associated with the null hypothesis that $\bX$ is not an outlier.   The sample complexity of estimation of $P(\bX\in \mathcal B)$, for all
$\mathcal B$, is equivalent to that of estimation of the density function $f(\bx)$.
Therefore, we adopt the mean integrated squared error (MISE) \cite{Tsybakov2009} as the
risk function: $\int E[(f_{\Omega}(\bx)-\hat{f}(\bx))^2]d\bx$ where $\hat{f}$ is an
empirical density estimator.  It is known that if $f$ is in the class of  Lipschitz
functions the minimax MISE risk is of the form \cite{Devroye&etal:96}: $\beta n^{-2/(1+p)}$, $\beta>0$. Hence
we use this minimax risk as a proxy for performance and the sample complexity for this
problem is: $n$ increases to infinity with $p$ at the asymptotic rate $\frac{p}{\log n}
=\alpha$, where $\alpha\in (0,\infty)$ is a constant. The critical threshold
(\ref{eq:rhocrit}) again converges to zero as $p\rightarrow \infty$.
\end{itemize}

 \begin{table}
\centering
{\small
\begin{tabular}{|c||c|c|c|c|c|}
   \hline
   {\bf {Task}}   & Screening & Detection & Support recovery  & Param. estimation & Perform. estimation\\
\hline \hline
   {\bf Risk} &$P(N_e>0)$& $P(N_e>0)$& $P({\mathrm {card}}\{{\calS} \Delta \hat{\calS}\}=\phi)$ & $E[\|\bOmega-\hat\bOmega\|^2_F]$ & $\int E[(f_{\Omega}(\bx)-\hat{f}(\bx))^2]d\bx$ \\
\hline
   {\bf Bound} &$1-e^{-\kappa_n}$& $pe^{-n\beta}$& $2^{p^\nu} e^{-n\beta}$ & $\frac{p\log p}{n}\beta$ & $n^{-2/(1+p)}\beta$\\
\hline
   {\bf Regimes} &$\frac{\log p}{n}\rightarrow \infty$& $\frac{\log p}{n}  \rightarrow \alpha$& $\frac{ p^\nu}{n}\rightarrow \alpha$ &  $\frac{p\log p}{n}\rightarrow \alpha$ & $\frac{p}{\log n}\rightarrow \alpha$  \\
\hline
   {\bf Threshold} &$\rho_c\rightarrow 1$& $\rho_c \rightarrow \rho^*$& $\rho_c\rightarrow 0$ & $\rho_c\rightarrow 0$  & $\rho_c\rightarrow 0$  \\
   \hline
\end{tabular}
}
\caption{\em \small Different sample complexity regimes characterize the difficulty of performing different inference tasks on the inverse covariance $\bOmega$. A similar table would hold for inference on the
covariance $\bSigma$. Task (1st row) specified risk function (2nd row) for which we can use the upper bound (3rd row).  The sample
complexity regimes (4th row) for which the bound remains constant depends on
the task: increasingly large sample sizes are required as the complexity of the task
increases from left to right.  The limiting value of the critical threshold $\rho_c$, defined in
(\ref{eq:rhocrit}), for screening (5th row) is shown for each regime. For the screening
task (detection of existence of large partial correlations) the bound is the purely high dimensional (large $p$ and fixed $n$) asymptotic limit of the probability of  false positives  associated with the test that at least one
variable is highly correlated (given by Theorem 1).   For the  detection task (detection
of existence of partial correlation of magnitude greater than $\rho \in (0,1)$) the bound is the mixed high dimensional (large $p$ and large $n$) asymptotic limit of the probability of false positives. In this regime, the critical threshold converges to $\rho^*$ where $0<\rho^*<\rho$.
For the support recovery task (model selection) the bound is the mixed high dimensional
asymptotic  probability of misclassification of the set of partial correlations using CONCORD. For this task the parameter $\nu$ is specified by {\em a priori} knowledge on sparsity of the support set: $\nu \in (0,1]$ and as the sparsity increases $\nu \rightarrow 0$.  For the estimation task the bound is the  mixed high dimensional asymptotic
squared Frobenius norm error on the MAP estimator of sparse inverse covariance. Finally, for the performance estimation task, the bound is the  mixed high dimensional asymptotic minmax bound on estimation MSE of the probability of any (Borel) uncertainty set. The values of the constants $\alpha$ and $\beta$ are not necessarily the same over different columns of the table.
}
\label{table:tasks}
\end{table}


We conclude this section with a comparison of computational complexity.  Unlike prediction and model selection, correlation screening methods are scalable to very high dimensions both in terms of computation and memory scalability. Popular sparse optimization approaches to covariance estimation and covariance selection are iterative and include penalized likelihood methods such as Glasso \cite{Glasso}, SPACE \cite{peng2009partial} and CONCORD  \cite{Khare2014} . The computational complexity of Glasso after $t$ iterations is of order $O(tp^3)$. This can be reduced when using regression based methods such as SPACE and CONCORD which have a computational complexity of order $\min \{O(tnp^2), O(tp^3) \}$.   In contrast, correlation and partial correlation screening are non-iterative algorithms and the computational complexity is only of order $O(np\log p)$ which can be considerably less than its penalized likelihood counterparts. The lower order  $O(np\log p)$ is due to the fact that building the thresholded sample covariance is equivalent to constructing a Euclidean ball graph over $p$ nodes in $n-1$ dimensional space, for which reliable approximate nearest neighbors (ANN) algorithms \cite{Arya&etal:ACM98} can be applied. Very fast and scalable C, Python and Matlab implementations of ANN algorithms are available, e.g., {\tt FLANN}, {\tt Gensim}, and {\tt annoy}, which have been implemented on datasets with $p$ in the millions and $n$ in the hundreds.

\section{Conclusions}
\label{sec:conclusions}
Big data is not just lots of data.  This monolithic characterization is overly simplistic and ignores the issues of inference,  limited samples, and reproducibility.
Big data is of limited utility without appropriate inferential tools, e.g., use of the dataset to produce empirical estimates, classifications, or decisions on the population that generated the data. Inferences in turn lack credibility without accounting for errors due to limited samples. Without credible inferences there is no reproducibility: another random sample from the sample population may produce completely different results.

This paper adopted a statistical perspective in which a large scale dataset is a set of $n$ random samples drawn from a population of $p$ variables and $p$ is large. We focused on the problem of correlation mining where the objective is to infer properties of the population covariance matrix from the samples. The reliability of the inferences from limited samples can be mathematically characterized by the high dimensional learning rates and sample complexity associated with the inference problem. These specify the relative rate at which $n$ must go to infinity as a function of $p$ in order to ensure accurate performance. The sample complexity falls into different high dimensional regimes including the classical regime, where $p$ is fixed and $n$ goes to infinity, the mixed dimensional where both $n$ and $p$ goes to infinity, and the purely high dimensional where $n$ is fixed and $p$ goes to infinity.


The comparative sampling complexity analysis illustrated in this paper unifies the problem of sample sizing for large scale inference problems and, in particular, for correlation mining. Indeed, different sample complexity regimes each occupy a niche for different correlation mining tasks. In particular, screening for high correlations is governed by purely high dimensional rates while model selection, covariance estimation, and uncertainty quantification require $n$ to go to infinity at progressively larger rates as a function of $p$. This implies, for example, that one can do screening with many fewer samples than are required for the other tasks. Furthermore, in situations where samples are acquired sequentially our analysis suggests that that one can adapt the inference task over time: starting with correlation screening when samples are few, and progressing on through support detection, covariance estimation, and uncertainty quantification as more and more samples are acquired. Such a strategy is explored in the context of the sequential prediction and regression criterion (SPARC) \cite{Firouzi&etal:AISTATS13}.

\section*{Acknowledgements}
The work of Alfred Hero was partially supported by US Air Force Office of Scientific Research grant award number FA9550-13-1-0043, US Army Research Office grant awards W911NF-11-1-0391 and W911NF-12-1-0443, US National Science Foundation award CCF-1217880, National Institutes of Health grant 2P01CA087634-06A2,  and the Consortium for Verification Technology under the US Department of Energy National Nuclear Security Administration, award DE-NA0002534. The work of Bala Rajaratnam was partially supported by US Air Force Office of Scientific Research grant award FA9550-13-1-0043, US National Science Foundation under grant DMS-0906392, DMS-CMG 1025465, AGS-1003823, DMS-1106642, DMS-CAREER-1352656, Defense Advanced Research Projects Agency DARPA-YFAN66001-111-4131, the UPS fund and SMC-DBNKY.

\bibliographystyle{sieee}
\bibliography{lib,libraryextra}
\newpage
\listoffigures
\end{document}

%% file: header_defns.tex
\newcommand{\be}{\begin{eqnarray}} \newcommand{\ee}{\end{eqnarray}}
\newcommand{\ben}{\begin{eqnarray*}} \newcommand{\een}{\end{eqnarray*}}
\newcommand{\bc}{\begin{center}} \newcommand{\ec}{\end{center}}
\newcommand{\bt}{\begin{tabbing}} \newcommand{\et}{\end{tabbing}}

\newcommand{\Reals}{\mbox{\rm I\kern-.25em R}}  
\newcommand{\Integers}{\mbox{\rm Z\kern-.25em Z}}  

\newcommand{\cov}{{\mathrm{cov}}}

\def\trace{{\mathrm{tr}}}

\def\kron{\bigotimes}

\def\diag{{\mathrm{diag}}}

\def\bS{{\mathbf{S}}}

\def\bP{{\mathbf{P}}}
\def\bA{{\mathbf{A}}}
\def\bB{{\mathbf{B}}}

\def\bD{{\mathbf{D}}}

\def\bI{{\mathbf{I}}}

\def\bP{{\mathbf{P}}}

\def\bR{{\mathbf{R}}}
\def\bS{{\mathbf{S}}}

\def\bW{{\mathbf{W}}}
\def\bx{{\mathbf{x}}}
\def\bX{{\mathbf{X}}}

\def\calS{{\mathcal{S}}}

\def\calS{{\mathcal{S}}}

\def\bR{{\mathbf{R}}}

\def\infinity{\infty}
\def\rank{{\mathrm{rank}}}
\def\logdet{{\mathrm {log det}}}
\def\bSigma{{\mathbf \Sigma}}

\def\bOmega{{\mathbf \Omega}}
\def\bmu{{\boldsymbol \mu}}

\def\calS{{\mathcal{S}}}

\def\calS{{\mathcal{S}}}

%% file: HR_ieee_proc14_submitted_final_arxiv_revised.bbl
\begin{thebibliography}{100}

\bibitem{Allen2010}
G.~I. Allen and R.~Tibshirani, ``Transposable regularized covariance models
  with an application to missing data imputation,'' {\em The Annals of Applied
  Statistics}, vol. 4, no. 2, pp. 764--790, 2010.

\bibitem{almasy1998multipoint}
L.~Almasy and J.~Blangero, ``Multipoint quantitative-trait linkage analysis in
  general pedigrees,'' {\em The American Journal of Human Genetics}, vol. 62,
  no. 5, pp. 1198--1211, 1998.

\bibitem{Anderson:03}
T.~W. Anderson, {\em An Introduction to Multivariate Statistical Analysis},
  Wiley, New York, 2003.

\bibitem{arratia1989tutorial}
R.~Arratia and L.~Gordon, ``Tutorial on large deviations for the binomial
  distribution,'' {\em Bulletin of mathematical biology}, vol. 51, no. 1, pp.
  125--131, 1989.

\bibitem{Arya&etal:ACM98}
S.~Arya, D.~M. Mount, N.~S. Netanyahu, R.~Silverman, and A.~Y. Wu, ``An optimal
  algorithm for approximate nearest neighbor searching fixed dimensions,'' {\em
  Journal of the ACM}, vol. 45, no. 6, pp. 891--923, 1998.

\bibitem{Bahadur1967}
R.~Bahadur, ``Rates of convergence of estimates and test statistics,'' {\em
  Annals of Mathematical Statistics}, vol. 38, pp. 303--324, 1967.

\bibitem{banerjee2006convex}
O.~Banerjee, L.~El~Ghaoui, A.~d'Aspremont, and G.~Natsoulis, ``Convex
  optimization techniques for fitting sparse {Gaussian} graphical models,'' in
  {\em ACM International Conference Proceeding Series}, volume 148, pp.
  89--96. Citeseer, 2006.

\bibitem{Banerjee2008a}
O.~Banerjee, L.~{El Ghaoui}, and A.~d'Aspremont, ``{Model Selection Through
  Sparse Maximum Likelihood Estimation for Multivariate Gaussian or Binary
  Data},'' {\em The Journal of Machine Learning Research}, vol. 9, pp.
  485--516, June 2008.

\bibitem{bartlett1998sample}
P.~L. Bartlett, ``The sample complexity of pattern classification with neural
  networks: the size of the weights is more important than the size of the
  network,'' {\em Information Theory, IEEE Transactions on}, vol. 44, no. 2,
  pp. 525--536, 1998.

\bibitem{BickelRitovTsybakov2009}
P.~Bickel, Y.~Ritov, and A.~Tsybakov, ``Simultaneous analysis of {Lasso} and
  {Dantzig} selector,'' {\em Annals of Statistics}, vol. 37, pp. 1705--1732,
  2009.

\bibitem{bickel&levina:aos08_2}
P.~Bickel and E.~Levina, ``Covariance regularization via thresholding,'' {\em
  Annals of Statistics}, vol. 34, no. 6, pp. 2577--2604, 2008.

\bibitem{biglieri2007mimo}
E.~Biglieri, R.~Calderbank, A.~Constantinides, A.~Goldsmith, A.~Paulraj, and
  H.~V. Poor, {\em MIMO wireless communications}, Cambridge University Press,
  2007.

\bibitem{bliss2002eim}
D.~Bliss, K.~Forsythe, A.~Hero~III, and A.~Yegulalp, ``{Environmental issues
  for MIMO capacity},'' {\em Signal Processing, IEEE Transactions on [see also
  Acoustics, Speech, and Signal Processing, IEEE Transactions on]}, vol. 50,
  no. 9, pp. 2128--2142, 2002.

\bibitem{bollen1998structural}
K.~A. Bollen, {\em Structural equation models}, Wiley Online Library, 1998.

\bibitem{Buhlmann:2011}
P.~B{\"u}hlmann and S.~van~de Geer, {\em Statistics for High-Dimensional Data:
  Methods, Theory and Applications}, Springer, 2011.

\bibitem{CandesTao2007}
E.~Cand\`{e}s and T.~Tao, ``The {Dantzig} selector: statistical estimation when
  p is much larger than n,'' {\em Annals of Statistics}, vol. 35, pp.
  2313--2351, 2007.

\bibitem{chamberland2003decentralized}
J.-F. Chamberland and V.~V. Veeravalli, ``Decentralized detection in sensor
  networks,'' {\em Signal Processing, IEEE Transactions on}, vol. 51, no. 2,
  pp. 407--416, 2003.

\bibitem{chellappa1985classification}
R.~Chellappa and S.~Chatterjee, ``Classification of textures using gaussian
  markov random fields,'' {\em Acoustics, Speech and Signal Processing, IEEE
  Transactions on}, vol. 33, no. 4, pp. 959--963, 1985.

\bibitem{chen2011robust}
Y.~Chen, A.~Wiesel, and A.~O. Hero, ``Robust shrinkage estimation of
  high-dimensional covariance matrices,'' {\em Signal Processing, IEEE
  Transactions on}, vol. 59, no. 9, pp. 4097--4107, 2011.

\bibitem{Chernoff1956}
H.~Chernoff, ``Large-sample theory: Parametric case,'' {\em Annals of
  Mathematical Statistics}, vol. 27, pp. 1--22, 1956.

\bibitem{chung2007detection}
P.-J. Chung, J.~F. B{\"o}hme, C.~F. Mecklenbrauker, and A.~O. Hero, ``Detection
  of the number of signals using the benjamini-hochberg procedure,'' {\em
  Signal Processing, IEEE Transactions on}, vol. 55, no. 6, pp. 2497--2508,
  2007.

\bibitem{craddock2013imaging}
R.~C. Craddock, S.~Jbabdi, C.-G. Yan, J.~T. Vogelstein, F.~X. Castellanos,
  A.~Di~Martino, C.~Kelly, K.~Heberlein, S.~Colcombe, and M.~P. Milham,
  ``Imaging human connectomes at the macroscale,'' {\em Nature methods}, vol.
  10, no. 6, pp. 524--539, 2013.

\bibitem{Cramer1946b}
H.~Cram\'{e}r, ``A contribution to the theory of statistical estimation,'' {\em
  Scandinavian Actuarial Journal}, vol. 29, pp. 85--94, 1946.

\bibitem{Cramer1946a}
H.~Cram\'{e}r, {\em Mathematical Methods of Statistics}, Princeton University
  Press, Princeton, NJ, 1946.

\bibitem{PhysRevLett.63.105}
J.~P. Crutchfield and K.~Young, ``Inferring statistical complexity,'' {\em
  Phys. Rev. Lett.}, vol. 63, pp. 105--108, Jul 1989.

\bibitem{dal2014making}
R.~Dal-R{\'e}, J.~P. Ioannidis, M.~B. Bracken, P.~A. Buffler, A.-W. Chan, E.~L.
  Franco, C.~La~Vecchia, and E.~Weiderpass, ``Making prospective registration
  of observational research a reality,'' {\em Science translational medicine},
  vol. 6, no. 224, pp. 224cm1--224cm1, 2014.

\bibitem{DalalRaj2014}
O.~Dalal and B.~Rajaratnam, ``{G-AMA: Sparse Gaussian Graphical Model
  Estimation via Alternating Minimization},'' {\em Technical Report, Department
  of Statistics, Stanford University (in revision)}, 2014.

\bibitem{dasgupta2005coarse}
S.~Dasgupta, ``Coarse sample complexity bounds for active learning,'' in {\em
  Advances in neural information processing systems}, pp.  235--242, 2005.

\bibitem{dawid_lauritzen}
A.~Dawid and S.~Lauritzen, ``{Hyper Markov laws in the statistical analysis of
  decomposable graphical models.},'' {\em Ann. Stat.}, vol. 21, no. 3, pp.
  1272--1317, 1993.

\bibitem{de2004discovery}
A.~De~La~Fuente, N.~Bing, I.~Hoeschele, and P.~Mendes, ``Discovery of
  meaningful associations in genomic data using partial correlation
  coefficients,'' {\em Bioinformatics}, vol. 20, no. 18, pp. 3565--3574, 2004.

\bibitem{Dempster1972}
A.~P. Dempster, ``{Covariance Selection},'' {\em Biometrics}, vol. 28, no. 1,
  pp. 157--175, March 1972.

\bibitem{Devroye&etal:96}
L.~Devroye, L.~Gy{\"o}rfi, and G.~Lugosi, {\em A probabilistic theory of
  pattern recognition}, Springer-Verlag, New York NY, 1996.

\bibitem{Donoho2006}
D.~Donoho, ``For most large underdetermined systems of linear equations the
  minimal $\ell_1$-norm solution is also the sparsest solution,'' {\em
  Communications on Pure and Applied Mathematics}, vol. 59, pp. 797--829, 2006.

\bibitem{Dutilleul1999}
P.~Dutilleul, ``The mle algorithm for the matrix normal distribution,'' {\em
  Journal of Statistical Computation and Simulation}, vol. 64, pp. 105--123,
  1999.

\bibitem{Efron1982}
B.~Efron, ``Maximum likelihood and decision theory,'' {\em Annals of
  Statistics}, vol. 10, pp. 340--356, 1982.

\bibitem{fan2008sure}
J.~Fan and J.~Lv, ``Sure independence screening for ultrahigh dimensional
  feature space,'' {\em Journal of the Royal Statistical Society: Series B
  (Statistical Methodology)}, vol. 70, no. 5, pp. 849--911, 2008.

\bibitem{Firouzi&etal:AISTATS13}
H.~Firouzi, A.~Hero, and B.~Rajaratnam, ``Predictive correlation screening:
  Application to two-stage predictor design in high dimension,'' in {\em
  Proceedings of AISTATS. Also available as arxiv:1303.2378}, 2013.

\bibitem{FirouziHeroRajaratnam2014}
H.~Firouzi, A.~Hero, and B.~Rajaratnam, ``Two-stage sampling, prediction and
  adaptive regression via correlation screening (sparcs),'' {\em arxiv
  1502:06189}, 2015.

\bibitem{Firouzi&etal:GlobalSIP13}
H.~Firouzi, D.~Wei, and A.~Hero, ``Spatio-temporal analysis of gaussian wss
  processes via complex correlation and partial correlation screening,'' in
  {\em Proceedings of IEEE GlobalSIP Conference. Also available as
  arxiv:1303.2378}, 2013.

\bibitem{firouziExcursions2014}
H.~Firouzi, D.~Wei, and A.~Hero, ``Spectral correlation hub screening of
  multivariate time series,'' in {\em Excursions in Harmonic Analysis: The
  February Fourier Talks at the Norbert Wiener Center}, R.~Balan, M.~Begu\'{e},
  J.~J. Benedetto, W.~Czaja, and K.~Okoudjou, editors, Springer, 2014.

\bibitem{firouzi2013local}
H.~Firouzi and A.~O. Hero, ``Local hub screening in sparse correlation
  graphs,'' in {\em SPIE Optical Engineering+ Applications}, pp.
  88581H--88581H. International Society for Optics and Photonics, 2013.

\bibitem{Fisher1922}
R.~Fisher, ``On the mathematical foundations of theoretical statistics,'' {\em
  Philosophical Transactions of the Royal Society of London, Series A}, vol.
  222, pp. 309--368, 1922.

\bibitem{Fisher1925}
R.~Fisher, ``Theory of statistical estimation,'' {\em Proceedings of the
  Cambridge Philosophical Society}, vol. 22, pp. 700--725, 1925.

\bibitem{fornell1982two}
C.~Fornell and F.~L. Bookstein, ``Two structural equation models: Lisrel and
  pls applied to consumer exit-voice theory,'' {\em Journal of Marketing
  research}, pp.  440--452, 1982.

\bibitem{Glasso}
J.~Friedman, T.~Hastie, and R.~Tibshirani, ``Sparse inverse covariance
  estimation with the graphical lasso,'' {\em Biostatistics}, vol. 9, no. 3,
  pp. 432--441, 2008.

\bibitem{Friedman2010}
J.~Friedman, T.~Hastie, and R.~Tibshirani.
\newblock {\em {Applications of the lasso and grouped lasso to the estimation
  of sparse graphical models}}, 2010.

\bibitem{Friedman&etal:08}
J.~Friedman, R.~Hastie, and R.~Tibshirani, ``Sparse inverse covariance
  estimation with the graphical lasso,'' {\em Biostatistics}, vol. 9, no. 3,
  pp. 432--441, 2008.

\bibitem{fuhrmann2008transmit}
D.~R. Fuhrmann and G.~San~Antonio, ``Transmit beamforming for mimo radar
  systems using signal cross-correlation,'' {\em Aerospace and Electronic
  Systems, IEEE Transactions on}, vol. 44, no. 1, pp. 171--186, 2008.

\bibitem{geman&geman:1984}
S.~Geman and D.~Geman, ``Stochastic relaxation, gibbs distributions, and the
  bayesian restoration of images,'' {\em Pattern Analysis and Machine
  Intelligence, IEEE Transactions on}, no. 6, pp. 721--741, 1984.

\bibitem{ghaoui2003worst}
L.~E. Ghaoui, M.~Oks, and F.~Oustry, ``Worst-case value-at-risk and robust
  portfolio optimization: A conic programming approach,'' {\em Operations
  Research}, vol. 51, no. 4, pp. 543--556, 2003.

\bibitem{gini2002covariance}
F.~Gini and M.~Greco, ``Covariance matrix estimation for cfar detection in
  correlated heavy tailed clutter,'' {\em Signal Processing}, vol. 82, no. 12,
  pp. 1847--1859, 2002.

\bibitem{Greenewald2013}
K.~Greenewald, T.~Tsiligkaridis, and A.~O. Hero, ``Kronecker sum decompositions
  of space-time data,'' in {\em Computational Advances in Multi-Sensor Adaptive
  Processing (CAMSAP), 2013 IEEE 5th International Workshop on}, pp.  65--68.
  IEEE, 2013.

\bibitem{guerci2003space}
J.~R. Guerci, {\em Space-time adaptive processing for radar}, Artech House,
  2003.

\bibitem{guerci2000optimal}
J.~Guerci, J.~Goldstein, and I.~Reed, ``Optimal and adaptive reduced-rank
  stap,'' {\em Aerospace and Electronic Systems, IEEE Transactions on}, vol.
  36, no. 2, pp. 647--663, 2000.

\bibitem{GRE2014}
D.~Guillot, B.~Rajaratnam, and J.~Emile-Geay, ``Statistical paleoclimate
  reconstructions via markov random fields,'' {\em Annals of Applied Statistics
  (to appear)}, 2014.

\bibitem{Guillot2012}
D.~Guillot, B.~Rajaratnam, B.~T. Rolfs, A.~Maleki, and I.~Wong, ``{Iterative
  Thresholding Algorithm for Sparse Inverse Covariance Estimation},'' in {\em
  Advances in Neural Information Processing Systems 25}, 2012.

\bibitem{haussler1994bounds}
D.~Haussler, M.~Kearns, and R.~E. Schapire, ``Bounds on the sample complexity
  of bayesian learning using information theory and the vc dimension,'' {\em
  Machine learning}, vol. 14, no. 1, pp. 83--113, 1994.

\bibitem{hero2012hub}
A.~Hero and B.~Rajaratnam, ``Hub discovery in partial correlation models,''
  {\em IEEE Trans. on Inform. Theory}, vol. 58, no. 9, pp. 6064--6078, 2012.
\newblock available as Arxiv preprint arXiv:1109.6846.

\bibitem{Hero:NIPS06}
A.~O. Hero, ``Geometric entropy minimization {(GEM)} for anomaly detection and
  localization,'' in {\em Proc. Neural Information Processing Systems (NIPS)
  Conference}, 2006.

\bibitem{Hero&Delap:Haykin95}
A.~O. Hero and R.~Delap, ``Task specific criteria for adaptive beamsumming with
  slow fading signals,'' in {\em Advances in Spectrum Analysis and Array
  Processing: Vol. III}, S.~Haykin, editor, Prentice Hall, Englewood-Cliffs,
  NJ, 1995.

\bibitem{hero2011large}
A.~Hero and B.~Rajaratnam, ``Large-scale correlation screening,'' {\em Journal
  of the American Statistical Association}, vol. 106, no. 496, pp. 1540--1552,
  2011.

\bibitem{hero2003secure}
A.~O. Hero, ``Secure space-time communication,'' {\em Information Theory, IEEE
  Transactions on}, vol. 49, no. 12, pp. 3235--3249, 2003.

\bibitem{Hero&Rajaratnam:Cambridge14}
A.~Hero and B.~Rajaratnam, ``Large scale correlation mining for biomolecular
  network discovery,'' in {\em Big data over networks}, S.~Cui, A.~Hero,
  Z.~Luo, and J.~Moura, editors. Cambridge Univ Press, 2015. Preprint available
  in Stanford University Dept. of Statistics Report series.

\bibitem{hsiao2014social}
K.-J. Hsiao, A.~Kulesza, and A.~Hero, ``Social collaborative retrieval,'' in
  {\em Proceedings of the 7th ACM international conference on Web search and
  data mining}, pp.  293--302. ACM, 2014.

\bibitem{Hsieh2011}
C.-J. Hsieh, M.~A. Sustik, I.~S. Dhillon, and P.~K. Ravikumar, ``{Sparse
  Inverse Covariance Matrix Estimation Using Quadratic Approximation},'' in
  {\em Advances in Neural Information Processing Systems 24}, 2011.

\bibitem{Johnson&Dudgeon:93}
D.~H. Johnson and D.~E. Dudgeon, {\em Array Signal Processing}, Prentice Hall,
  Englewood-Cliffs N.J., 1993.

\bibitem{kakade2003sample}
S.~M. Kakade et~al., {\em On the sample complexity of reinforcement learning},
  PhD thesis, University of London, 2003.

\bibitem{elkaroui:annstat08}
N.~E. Karoui, ``Operator norm consistent estimation of large dimensional sparse
  covariance matrices,'' {\em Annals of Statistics}, vol. to appear, , 2008.

\bibitem{Kay:91}
S.~M. Kay, {\em Statistical Estimation}, Prentice-Hall, Englewood-Cliffs N.J.,
  1991.

\bibitem{Kay:98}
S.~Kay, {\em Fundamentals of Statistical Signal Processing, Volume 2: Detection
  Theory}, Prentice-Hall, Englewood-Cliffs N.J., 1998.

\bibitem{Kelly&Forsythe:TR848}
E.~J. Kelly and K.~M. Forsythe, ``Adaptive detection and parameter estimation
  for multidimensional signal models,'' Technical Report 848, M.I.T. Lincoln
  Laboratory, April, 1989.

\bibitem{KhareOhRaj2014}
K.~Khare, S.~Oh, and B.~Rajaratnam, ``{A convex pseudo-likelihood framework for
  high dimensional partial correlation estimation with convergence
  guarantees},'' {\em Journal of the Royal Statistical Society: Series B
  (Statistical Methodology), to appear}, 2014.

\bibitem{Khare2014}
K.~Khare, S.~Oh, and B.~Rajaratnam, ``{A convex pseudo-likelihood framework for
  high dimensional partial correlation estimation with convergence
  guarantees},'' {\em Journal of the Royal Statistical Society: Series B
  (Statistical Methodology), to appear}, 2014.

\bibitem{KRopt2014}
K.~Khare and B.~Rajaratnam, ``{ Convergence of cyclic coordinatewise l1
  minimization},'' Technical report, Stanford University, 2014.

\bibitem{Khare2011}
K.~Khare and B.~Rajaratnam, ``{Wishart distributions for decomposable
  covariance graph models},'' {\em The Annals of Statistics}, vol. 39, no. 1,
  pp. 514--555, March 2011.

\bibitem{KieferWolfowitz1956}
J.~Kiefer and J.~Wolfowitz, ``Consistency of the maximum likelihood esitmator
  in the presence of infinitely many incidental parameters,'' {\em Annals of
  Mathematical Statistics}, vol. 27, pp. 887--906, 1956.

\bibitem{kim2001comparison}
H.~S. Kim and A.~O. Hero, ``Comparison of glr and invariant detectors under
  structured clutter covariance,'' {\em Image Processing, IEEE Transactions
  on}, vol. 10, no. 10, pp. 1509--1520, 2001.

\bibitem{koren2009matrix}
Y.~Koren, R.~Bell, and C.~Volinsky, ``Matrix factorization techniques for
  recommender systems,'' {\em Computer}, vol. 42, no. 8, pp. 30--37, 2009.

\bibitem{Korrat:2007ec}
A.~Korrat, T.~Greiner, M.~Maurer, T.~Metz, and H.-H. Fiebig, ``Gene
  signature-based prediction of tumor response to cyclophosphamide.,'' {\em
  Cancer Genomics Proteomics}, vol. 4, no. 3, pp. 187--195, 2007.

\bibitem{labrinidis2012challenges}
A.~Labrinidis and H.~Jagadish, ``Challenges and opportunities with big data,''
  {\em Proceedings of the VLDB Endowment}, vol. 5, no. 12, pp. 2032--2033,
  2012.

\bibitem{lauritzen1996graphical}
S.~Lauritzen, {\em Graphical models}, volume~17, Oxford University Press, USA,
  1996.

\bibitem{LeCam1953}
L.~Le~Cam, ``On some asymptotic properties of maximum likelihood estimates and
  related {Bayes'} estimates,'' {\em University of California publications in
  statistics}, vol. 1, pp. 277--330, 1953.

\bibitem{LeCam1986}
L.~Le~Cam, {\em Asymptotic Methods in Statistical Decision Theory},
  Springer-Verlag, New York, 1986.

\bibitem{ledoit2003improved}
O.~Ledoit and M.~Wolf, ``Improved estimation of the covariance matrix of stock
  returns with an application to portfolio selection,'' {\em Journal of
  empirical finance}, vol. 10, no. 5, pp. 603--621, 2003.

\bibitem{ledoit2004honey}
O.~Ledoit and M.~Wolf, ``Honey, i shrunk the sample covariance matrix,'' {\em
  The Journal of Portfolio Management}, vol. 30, no. 4, pp. 110--119, 2004.

\bibitem{Ledoit&Wolf:JMA04}
Q.~Ledoit and M.~Wolf, ``A well conditioned estimator for large dimensional co-
  variance matrices,'' {\em J. Multiv. Anal.}, vol. 88, pp. 365--411, 2004.

\bibitem{LeeHastie2014}
J.~Lee and T.~Hastie, ``{Learning the structure of mixed graphical models},''
  {\em Journal of Computational and Graphical Statistics, to appear)}, 2014.

\bibitem{letac_massam_2007}
G.~Letac and H.~Massam, ``{Wishart distributions for decomposable graphs},''
  {\em Annals of Statistics}, vol. 35, no. 3, , 2007.

\bibitem{li2009mimo}
J.~Li and P.~Stoica, {\em MIMO Radar Signal Processing}, John Wiley \& Sons,
  Inc., Hoboken, NJ, 2009.

\bibitem{lynch2008big}
C.~Lynch, ``Big data: How do your data grow?,'' {\em Nature}, vol. 455, no.
  7209, pp. 28--29, 2008.

\bibitem{madigan2013evaluating}
D.~Madigan, P.~B. Ryan, M.~Schuemie, P.~E. Stang, J.~M. Overhage, A.~G.
  Hartzema, M.~A. Suchard, W.~DuMouchel, and J.~A. Berlin, ``Evaluating the
  impact of database heterogeneity on observational study results,'' {\em
  American journal of epidemiology}, vol. 178, no. 4, pp. 645--651, 2013.

\bibitem{mardia1988multi}
K.~Mardia, ``Multi-dimensional multivariate gaussian markov random fields with
  application to image processing,'' {\em Journal of Multivariate Analysis},
  vol. 24, no. 2, pp. 265--284, 1988.

\bibitem{marx2013biology}
V.~Marx, ``Biology: The big challenges of big data,'' {\em Nature}, vol. 498,
  no. 7453, pp. 255--260, 2013.

\bibitem{mcintosh1996spatial}
A.~McIntosh, F.~Bookstein, J.~V. Haxby, and C.~Grady, ``Spatial pattern
  analysis of functional brain images using partial least squares,'' {\em
  Neuroimage}, vol. 3, no. 3, pp. 143--157, 1996.

\bibitem{Meinshausen2006a}
N.~Meinshausen and P.~B\"{u}hlmann, ``{High-dimensional graphs and variable
  selection with the Lasso},'' {\em The Annals of Statistics}, vol. 34, no. 3,
  pp. 1436--1462, 2006.

\bibitem{michener2012ecoinformatics}
W.~K. Michener and M.~B. Jones, ``Ecoinformatics: supporting ecology as a
  data-intensive science,'' {\em Trends in ecology \& evolution}, vol. 27, no.
  2, pp. 85--93, 2012.

\bibitem{Morrison:90}
D.~F. Morrison, {\em Multivariate statistical methods}, McGraw Hill, New York,
  1990.

\bibitem{Muirhead:82}
R.~J. Muirhead, {\em Aspects of Multivariate Statistical Theory}, Wiley, New
  York, 1982.

\bibitem{NeymanPearson1933}
J.~Neyman and E.~Pearson, ``On the problem of the most efficient tests of
  statistical hypotheses,'' {\em Philosophical Transactions of the Royal
  Society of London, Series A}, vol. 231, pp. 289--337, 1933.

\bibitem{niyogi1996relationship}
P.~Niyogi and F.~Girosi, ``On the relationship between generalization error,
  hypothesis complexity, and sample complexity for radial basis functions,''
  {\em Neural Computation}, vol. 8, no. 4, pp. 819--842, 1996.

\bibitem{ODKRNips14}
S.~Oh, O.~Dalal, K.~Khare, and B.~Rajaratnam, ``{Optimization methods for
  sparse pseudo-likelihood graphical model selection},'' in {\em Advances in
  Neural Information Processing Systems 27}, 2014.

\bibitem{patwari2005locating}
N.~Patwari, J.~N. Ash, S.~Kyperountas, A.~O. Hero, R.~L. Moses, and N.~S.
  Correal, ``Locating the nodes: cooperative localization in wireless sensor
  networks,'' {\em Signal Processing Magazine, IEEE}, vol. 22, no. 4, pp.
  54--69, 2005.

\bibitem{patwari2003relative}
N.~Patwari, A.~O. Hero, M.~Perkins, N.~S. Correal, and R.~J. O'dea, ``Relative
  location estimation in wireless sensor networks,'' {\em Signal Processing,
  IEEE Transactions on}, vol. 51, no. 8, pp. 2137--2148, 2003.

\bibitem{peng2009a}
J.~Peng, P.~Wang, N.~Zhou, and J.~Zhu, ``{Partial Correlation Estimation by
  Joint Sparse Regression Models},'' {\em Journal of the American Statistical
  Association}, vol. 104, no. 486, pp. 735--746, 2009.

\bibitem{peng2009partial}
J.~Peng, P.~Wang, N.~Zhou, and J.~Zhu, ``Partial correlation estimation by
  joint sparse regression models,'' {\em Journal of the American Statistical
  Association}, vol. 104, no. 486, , 2009.

\bibitem{rajapakse2010networking}
I.~Rajapakse, D.~Scalzo, S.~J. Tapscott, S.~T. Kosak, and M.~Groudine,
  ``Networking the nucleus,'' {\em Molecular systems biology}, vol. 6, no. 1, ,
  2010.

\bibitem{Rajaratnam2008}
B.~Rajaratnam, H.~Massam, and C.~M. Carvalho, ``{Flexible covariance estimation
  in graphical Gaussian models},'' {\em The Annals of Statistics}, vol. 36, no.
  6, pp. 2818--2849, December 2008.

\bibitem{Rao1947}
C.~Rao, ``Large sample tests of statistical hypotheses concerning several
  parameters with applications to problems of estimation,'' {\em Mathematical
  Proceedings of the Cambridge Philosophical Society}, vol. 44, pp. 50--57,
  1947.

\bibitem{Rao1963}
C.~Rao, ``Criteria of estimation in large samples,'' {\em Sankhy\={a}: The
  Indian Journal of Statistics, Series A}, vol. 25, pp. 189--206, 1963.

\bibitem{Reed&Yu:ASSP90}
I.~S. Reed and X.~Yu, ``Adaptive multi-band {CFAR} detection of an optical
  pattern with unknown spectral distribution,'' {\em IEEE Trans. Acoust.,
  Speech, and Sig. Proc.}, vol. 38, no. 10, pp. 1760--1771, 1990.

\bibitem{rissanen1989stochastic}
J.~Rissanen, {\em Stochastic complexity in statistical inquiry theory}, World
  Scientific Publishing Co., Inc., 1989.

\bibitem{Robey&etal:AES92}
F.~Robey, D.~Fuhrmann, E.~Kelly, and R.~Nitzberg, ``A {CFAR} adaptive matched
  filter detector,'' {\em {IEEE} Transactions on Aerospace and Electronic
  Systems}, vol. 43, no. 12, pp. 2964--2974, l992.

\bibitem{Rocha2008}
G.~Rocha, P.~Zhao, and B.~Yu, ``{A path following algorithm for Sparse
  Pseudo-Likelihood Inverse Covariance Estimation (SPLICE)},'' Technical
  report, Statistics Department, UC Berkeley, Berkeley, CA, 2008.

\bibitem{Rothman2009}
A.~J. Rothman, E.~Levina, and J.~Zhu, ``{Generalized thresholding of Large
  Covariance Matrices},'' {\em Journal of the American Statistical
  Association}, vol. 104(485), pp. 177--186, 2009.

\bibitem{rothman2008sparse}
A.~Rothman, P.~Bickel, E.~Levina, and J.~Zhu, ``Sparse permutation invariant
  covariance estimation,'' {\em Electronic Journal of Statistics}, vol. 2, pp.
  494--515, 2008.

\bibitem{rudin2014discovery}
C.~Rudin, D.~Dunson, R.~Irizarry, H.~Ji, E.~Laber, J.~Leek, T.~McCormick,
  S.~Rose, C.~Schafer, M.~van~der Laan, et~al., ``Discovery with data:
  Leveraging statistics with computer science to transform science and
  society,'' {\em White Papre, American Statistical Association}, July 2014.

\bibitem{Schafer&etal:PhysicsofScales}
B.~M. Sch\"{a}fer and M.~Bertelman, {\em Physics of scales}, Zentrum f\"{u}r
  Astronomie, Universt\"{a}t Hiedelberg, 2013.

\bibitem{schafer2005shrinkage}
J.~Sch{\"a}fer and K.~Strimmer, ``A shrinkage approach to large-scale
  covariance matrix estimation and implications for functional genomics,'' {\em
  Statistical applications in genetics and molecular biology}, vol. 4, no. 1, ,
  2005.

\bibitem{Schneider01}
T.~{Schneider}, ``{Analysis of incomplete climate data: estimation of mean
  values and covariance matrices and imputation of missing values.},'' {\em J.
  Clim.}, vol. 14, pp. 853--871, 2001.

\bibitem{Schwartz:IT69}
R.~E. Schwartz, ``Minimax {CFAR} detection in additive {G}aussian noise of
  unknown covariance,'' {\em IEEE Trans. on Inform. Theory}, vol. 15, no. 4,
  pp. 722--725, July 1969.

\bibitem{Scott&Nowak:JMLR06}
C.~Scott and R.~Nowak, ``Learning minimum volume sets,'' {\em Journal of
  Machine Learning Research}, vol. 7, pp. 665--704, April 2006.

\bibitem{Smerdon_et_al2010}
J.~E. Smerdon, A.~Kaplan, D.~Chang, and M.~N. Evans, ``{A pseudoproxy
  evaluation of the CCA and RegEM methods for reconstructing climate fields of
  the last millennium},'' {\em Journal of Climate}, vol. 23, no. 18, pp.
  4856--4880, 2010.

\bibitem{sommerfeld1949partial}
A.~Sommerfeld, {\em Partial differential equations in physics}, volume~1,
  Academic press, 1949.

\bibitem{spiegelhalter2002bayesian}
D.~J. Spiegelhalter, N.~G. Best, B.~P. Carlin, and A.~Van Der~Linde, ``Bayesian
  measures of model complexity and fit,'' {\em Journal of the Royal Statistical
  Society: Series B (Statistical Methodology)}, vol. 64, no. 4, pp. 583--639,
  2002.

\bibitem{sporns2005human}
O.~Sporns, G.~Tononi, and R.~K{\"o}tter, ``The human connectome: a structural
  description of the human brain,'' {\em PLoS computational biology}, vol. 1,
  no. 4, pp. e42, 2005.

\bibitem{stoica2008using}
P.~Stoica, J.~Li, X.~Zhu, and J.~R. Guerci, ``On using a priori knowledge in
  space-time adaptive processing,'' {\em Signal Processing, IEEE Transactions
  on}, vol. 56, no. 6, pp. 2598--2602, 2008.

\bibitem{todros2012measure}
K.~Todros and A.~O. Hero, ``On measure transformed canonical correlation
  analysis,'' {\em Signal Processing, IEEE Transactions on}, vol. 60, no. 9,
  pp. 4570--4585, 2012.

\bibitem{VanTrees}
H.~L.~V. Trees, {\em Detection, Estimation and Modulation Theory: Part I}, John
  Wiley \& Sons, 2001.

\bibitem{trelles2011big}
O.~Trelles, P.~Prins, M.~Snir, and R.~C. Jansen, ``Big data, but are we
  ready?,'' {\em Nature Reviews Genetics}, vol. 12, no. 3, pp. 224--224, 2011.

\bibitem{Tsiligkaridis&Hero:TSP2013Kroneckerdecomposition}
T.~Tsiligkaridis and A.~Hero, ``Covariance estimation in high dimensions via
  kronecker product expansions,'' {\em IEEE Trans. on Signal Processing (also
  available as arXiv:1302.2686)}, vol. 61, no. 21, pp. 5347--5360, 2013.

\bibitem{tsiligkaridis2013convergenceTSP}
T.~Tsiligkaridis, A.~Hero, and S.~Zhou, ``Convergence properties of kronecker
  graphical lasso algorithms,'' {\em Signal Processing, IEEE Transactions on},
  vol. 61, no. 7, pp. 1743--1755, 2013.

\bibitem{Tsybakov2009}
A.~Tsybakov, {\em Introduction to nonparametric estimation}, Springer Verlag,
  2009.

\bibitem{van2010exploring}
M.~P. Van Den~Heuvel and H.~E. Hulshoff~Pol, ``Exploring the brain network: a
  review on resting-state fmri functional connectivity,'' {\em European
  Neuropsychopharmacology}, vol. 20, no. 8, pp. 519--534, 2010.

\bibitem{van2012human}
D.~C. Van~Essen, K.~Ugurbil, E.~Auerbach, D.~Barch, T.~Behrens, R.~Bucholz,
  A.~Chang, L.~Chen, M.~Corbetta, S.~W. Curtiss, et~al., ``The human connectome
  project: a data acquisition perspective,'' {\em Neuroimage}, vol. 62, no. 4,
  pp. 2222--2231, 2012.

\bibitem{van1993approximation}
C.~F. Van~Loan and N.~Pitsianis, {\em Approximation with Kronecker products},
  Springer, 1993.

\bibitem{van1997localization}
B.~D. Van~Veen, W.~van Drongelen, M.~Yuchtman, and A.~Suzuki, ``Localization of
  brain electrical activity via linearly constrained minimum variance spatial
  filtering,'' {\em Biomedical Engineering, IEEE Transactions on}, vol. 44, no.
  9, pp. 867--880, 1997.

\bibitem{vardi1996network}
Y.~Vardi, ``Network tomography: Estimating source-destination traffic
  intensities from link data,'' {\em Journal of the American Statistical
  Association}, vol. 91, no. 433, pp. 365--377, 1996.

\bibitem{Wainwright2009a}
M.~Wainwright, ``Information-theoretic limitations on sparsity recovery in the
  high-dimensional and noisy setting,'' {\em IEEE Transactions on Information
  Theory}, vol. 55, pp. 5728--5741, 2009.

\bibitem{Wainwright2009b}
M.~Wainwright, ``Sharp thresholds for high-dimensional and noisy sparsity
  recovery using $\ell_1$-constrained quadratic programming ({Lasso}),'' {\em
  IEEE Transactions on Information Theory}, vol. 55, pp. 2183--2202, 2009.

\bibitem{wainwright2008graphical}
M.~Wainwright and M.~Jordan, ``Graphical models, exponential families, and
  variational inference,'' {\em Foundations and Trends{\textregistered} in
  Machine Learning}, vol. 1, no. 1-2, pp. 1--305, 2008.

\bibitem{Wald1941a}
A.~Wald, ``Asymptotically most powerful tests of statistical hypotheses,'' {\em
  Annals of Mathematical Statistics}, vol. 12, pp. 1--19, 1941.

\bibitem{Wald1941b}
A.~Wald, ``Some examples of asymptotically most powerful tests,'' {\em Annals
  of Mathematical Statistics}, vol. 12, pp. 396--408, 1941.

\bibitem{Wald1943}
A.~Wald, ``Tests of statistical hypotheses concerning several parameters when
  the number of observations is large,'' {\em Transactions of the American
  Mathematical Society}, vol. 54, pp. 426--482, 1943.

\bibitem{Wald1949}
A.~Wald, ``Note on the consistency of the maximum likelihood estimate,'' {\em
  Annals of Mathematical Statistics}, vol. 20, pp. 595--601, 1949.

\bibitem{WEGSR2014}
J.~Wang, J.~Emile-Geay, D.~Guillot, J.~E. Smerdon, and B.~Rajaratnam,
  ``Evaluating climate field reconstruction techniques using improved
  emulations of real-world conditions,'' {\em Climates of the Past}, vol. 10,
  pp. 1--19, 2014.

\bibitem{Werner2007}
K.~Werner and M.~Jansson, ``Estimation of kronecker structured channel
  covariances using training data,'' in {\em Proceedings of EUSIPCO}, 2007.

\bibitem{werner2008estimation}
K.~Werner, M.~Jansson, and P.~Stoica, ``On estimation of covariance matrices
  with kronecker product structure,'' {\em Signal Processing, IEEE Transactions
  on}, vol. 56, no. 2, pp. 478--491, 2008.

\bibitem{west2003bayesian}
M.~West, ``Bayesian factor regression models in the large p, small n
  paradigm,'' {\em Bayesian statistics}, vol. 7, no. 2003, pp. 723--732, 2003.

\bibitem{wicks2006space}
M.~C. Wicks, M.~Rangaswamy, R.~Adve, and T.~Hale, ``Space-time adaptive
  processing: a knowledge-based perspective for airborne radar,'' {\em Signal
  Processing Magazine, IEEE}, vol. 23, no. 1, pp. 51--65, 2006.

\bibitem{Wiesel&Hero:SP09}
A.~Wiesel and A.~O. Hero, ``Decomposable principal components analysis,'' {\em
  IEEE Trans. on Signal Processing}, vol. 57, no. 11, pp. 4369--4378, nov 2009.

\bibitem{Wilks1938}
S.~Wilks, ``The large-sample distribution of the likelihood ratio for testing
  composite hypotheses,'' {\em Annals of Mathematical Statistics}, vol. 9, pp.
  60--62, 1938.

\bibitem{Willsky:IEEE02}
A.~Willsky, ``Multi-resolution {M}arkov models for signal and image
  processing,'' {\em Proc. of IEEE}, pp.  1396--1458, 2002.

\bibitem{Xu&etal:ICC2009}
K.~Xu, Y.~Chen, M.~Kliger, P.~Woolf, and A.~Hero, ``Revealing social networks
  of spammers through spectral clustering,'' in {\em IEEE Intl Conf on
  Communications (ICC)}, Dresden, June 2009.

\bibitem{young2011deming}
S.~S. Young and A.~Karr, ``Deming, data and observational studies,'' {\em
  Significance}, vol. 8, no. 3, pp. 116--120, 2011.

\bibitem{ZhaoYu2006}
P.~Zhao and B.~Yu, ``On model selection consistency of lasso,'' {\em Journal of
  Machine Learning Research}, vol. 7, pp. 2541--2563, 2006.

\bibitem{zhu2005high}
D.~Zhu, A.~O. Hero, Z.~S. Qin, and A.~Swaroop, ``High throughput screening of
  co-expressed gene pairs with controlled false discovery rate (fdr) and
  minimum acceptable strength (mas),'' {\em Journal of Computational Biology},
  vol. 12, no. 7, pp. 1029--1045, 2005.

\bibitem{zhu2007bayesian}
D.~Zhu and A.~O. Hero~III, ``Bayesian hierarchical model for large-scale
  covariance matrix estimation,'' {\em Journal of Computational Biology}, vol.
  14, no. 10, pp. 1311--1326, 2007.

\end{thebibliography}
